\documentclass[11pt,twoside]{article} 
\usepackage{amssymb,amsmath} 
\usepackage{verbatim} 
\def \P{{\hbox{\vrule width 0.6pt height 6.8pt depth -.2pt\kern-0.2pt P}}}

\def \Lip{\hbox{Lip}}

\def \R {\mathbb R} 
 
\def \N {\mathbb N} 
\def \e {\varepsilon} 
\def\tilde{\widetilde}

\newcommand\Z{{\mathbb{Z}}} 
 
\def\ddq{\dot\Delta_q} 
\def\d{\partial}

\def\cP{{\cal P}} 
\def\cQ{{\cal Q}}

\def\cC{{\cal C}}

\def\div{\, \mbox{div}\,  } 
 
\def\sgn{\, \mbox{sgn}\,  } 
\def\Supp{\, \mbox{Supp}\,  } 
\newenvironment{p}{ 
\begin{description}  
\item\textit{\textbf{Proof:}}~} 
{\hfill\rule{2.1mm}{2.1mm} 
\end{description}} 
\newcommand{\Sum}{\displaystyle \sum} 
 
\newcommand{\Int}{\displaystyle \int}

\newcommand{\Sup}{\displaystyle \sup}

\newcommand{\du}{\delta\! u} 
\newcommand{\dU}{\delta\! U} 
\newcommand{\dt}{\delta\!\theta} 
\newcommand{\dTheta}{\delta\!\Theta} 
\newcommand\dPi{\delta\!\Pi} 
 
\def \epsilon {\varepsilon} 
\newtheorem{defi}{Definition}
\newtheorem{lem}{Lemma}

\newtheorem{pro}{Proposition}
\newtheorem{theo}{Theorem}
\newtheorem{rem}{Remark}
\newcommand{\remark}{\smallbreak\noindent{\bf Remark:~}}
 
\begin{document} 
\title{Existence and uniqueness results
for the Boussinesq system with data in Lorentz spaces} 
\date{}

\author{Rapha\"el  Danchin\footnote{Universit\'e Paris-Est, LAMA UMR 8050, 61 av. 
du G\'en\'eral de Gaulle, 94010 Cr\'eteil Cedex, France. E-mail: danchin$@$univ-paris12.fr} and 
Marius Paicu\footnote{Universit\'e Paris-Sud, Laboratoire de Math\'ematiques,
B\^atiment 425,  91405 Orsay Cedex, France.}}
\maketitle

\begin{abstract} 
This paper is devoted to the study of the Cauchy problem 
for the Boussinesq system with partial viscosity 
in dimension $N\geq3.$
First we prove 
a global  existence result for  data
in  Lorentz spaces satisfying
a smallness condition which 
is at the scaling of the equations. 
Second, we get
a  uniqueness result in Besov spaces
with {\it negative} indices of regularity
(despite the fact that there is
no smoothing effect on the temperature).
The proof relies
on a priori estimates with loss of regularity 
 for the nonstationary
Stokes system with convection. 
As a corollary, we obtain 
a global existence and uniqueness result
for small data in Lorentz spaces.
\end{abstract} 

 
\section{Introduction and main results} 

The present paper is devoted to the mathematical study
of the so-called Boussinesq system with partial viscosity:
\begin{equation}\label{eq:boussinesq} 
\begin{cases} 
\partial_t\theta+u\cdot\nabla \theta=0\\ 
\partial_t u+u\cdot\nabla u-\nu \Delta u+\nabla\Pi=\theta\, e_N\\ 
\div u=0. 
\end{cases}  
\end{equation} 
Above, $\theta=\theta(t,x)$ and $\Pi=\Pi(t,x)$
are real valued functions,
and $u=u(t,x)$ is a time dependent 
vector field. 
We denote by $e_N$  the vertical 
unit vector of $\R^N$.
It is assumed that the space variable $x$ belongs
to $\R^N.$
We supplement 
the system with Cauchy conditions $(\theta_0,u_0)$
at time $t=0$ and address the question of
solvability for  $t\geq0$. 
\smallbreak
Boussinesq system arises in simplified models for geophysics
in which case $u$ stands for the velocity field
and the forcing term $\theta\,e_N$ is proportional  
either to the temperature, or to the salinity
or to the density (see \cite{Ped}).
 Here, we shall call $\theta$ the temperature.

Remark that the standard incompressible Navier-Stokes equations
arise as a particular case of \eqref{eq:boussinesq} (just take
$\theta\equiv0$). Hence, it is tempting to study whether the classical 
theory  for the Navier-Stokes equations
may be extended to those
more general fluids. 
As far as one is concerned with \emph{global} results, the main difficulty that
we have to face is that $\theta$ is transported
by the flow of $u.$ Hence it has
no time decay nor gain of smoothness whatsoever
and the standard approach for solving Navier-Stokes equations
with a (given) source term is bound to fail.
Nevertheless, in a recent paper (see \cite{DANCHINPAICU}), 
we have stated that system \eqref{eq:boussinesq} has
global finite energy weak solutions 
for any data $(\theta_0,u_0)$ in $L^2$ and
that uniqueness holds true in dimension two.
That result may be compared with Leray's theorem
\cite{Leray} for the Navier-Stokes equations.

In  \cite{DANCHINPAICU},  we have also obtained a global 
existence result (for small data) 
in the spirit of Fujita and Kato's theorem (see \cite{FK}
and \cite{K}) for the Navier-Stokes equations:
\begin{theo}\label{th:globalsmooth} Let $N\geq3.$
Let  $\theta_0\in L^{\frac N3}(\R^N)\cap\dot B^0_{N,1}(\R^N)$
 and  $u_0\in L^{N,\infty}(\R^N)\cap\dot B^{-1+\frac Np}_{p,1}(\R^N)$
 for some  $p\in[N,\infty].$
There exists a positive $c$ depending only on $N$
and such that if 
$$
\|u_0\|_{L^{N,\infty}}+\nu^{-1}\|\theta_0\|_{L^{\frac N3}}\leq c\nu
$$
then  Boussinesq system $\eqref{eq:boussinesq}$ admits a  unique global
solution   $$ (\theta,u,\nabla\Pi)\in{\cal C}(\R_+;\dot B^0_{N,1})\times
\Bigl({\cal C}(\R_+;\dot B^{\frac Np-1}_{p,1}) \cap
L^1_{loc}(\R_+;\dot B^{\frac Np+1}_{p,1})\Bigr)^N \times
\Bigl(L^1_{loc}(\R_+;\dot B^{\frac Np-1}_{p,1})\Bigr)^N.
$$
Besides, from dimension $N=4$ on, 
space  $L^{\frac N3}(\R^N)$ may be 
replaced by $L^{\frac N3,\infty}(\R^N).$
\end{theo}
The reader is referred to section \ref{s:tools}
for the definition of Lorentz spaces $L^{p,\infty}$
and Besov spaces $\dot B^s_{p,r}$ which have been used
above. Remind that we have the (continuous) inclusions
$$
\dot B^0_{N,1}\subset L^N\subset L^{N,\infty}.
$$
Theorem \ref{th:globalsmooth} may be interpreted in
terms of {\it scaling of Boussinesq system}.
Indeed, \eqref{eq:boussinesq} is obviously invariant
under the transform
$$
u(t,x)\longmapsto\lambda u(\lambda^2t,\lambda x)\quad\text{and}\quad
\theta(t,x)\longmapsto\lambda^3\theta(\lambda^2t,\lambda x)\quad\hbox{for all}
\ \ \lambda>0.
$$
Hence data  $(\theta_0,u_0)$ belong to a functional space 
$E$ which is \emph{at the scaling of} the system
if and only if the norm of $E$ is  
invariant by 
\begin{equation}\label{eq:scaling}
u_0(x)\longmapsto\lambda u_0(\lambda x)\quad\text{and}\quad
\theta_0(x)\longmapsto\lambda^3\theta_0(\lambda x).
\end{equation}
In dimension $N\geq3,$ the spaces
$L^{\frac N3}\times\bigl(L^{N,\infty}\bigr)^N$
and $L^{\frac N3,\infty}\times\bigl(L^{N,\infty}\bigr)^N$
satisfy \eqref{eq:scaling}.
Hence, Theorem \ref{th:globalsmooth} is a global well-posedness
result for suitably smooth data satisfying 
a scaling invariant smallness condition.
In fact, prescribing additional regularity 
in homogeneous Besov spaces ensures that the  velocity belongs to
the set $L^1_{loc}(\R_+;\Lip)$ of
locally integrable functions in $t$ 
with values in the set of Lipschitz vector fields. 
This property was needed for 
uniqueness. 
Note in passing that the additional condition
required for the initial velocity still
satisfies \eqref{eq:scaling}. 
\medbreak
The present paper aims at  
weakening as much as possible the Besov space 
 assumption. Ideally, we would like to
  consider \emph{general} (small) data $(\theta_0,u_0)$ in 
the scaling invariant space 
$L^{\frac N3,\infty}\times\bigl(L^{N,\infty}\bigr)^N.$

For proving existence, the strategy is the following.
We smooth out the data so that Theorem \ref{th:globalsmooth}
provides a sequence $(\theta^n,u^n)$ of global solutions
which is bounded in
$$
L^\infty\bigl(\R_+;L^{\frac
  N3,\infty}\times\bigl(L^{N,\infty}\bigr)^N\bigr).
$$
Obviously, those bounds are insufficient to pass
to the limit in the nonlinear terms.
However, the parabolicity of the second equation of \eqref{eq:boussinesq}
provides additional regularity 
so that there is some chance to pass
to the limit anyway
provided $\theta_0$ belongs also to $L^{p,\infty}$
for some large enough $p.$
Those basic considerations will enable us to prove
the following statement:
\begin{theo}\label{th:existence}  
Let $u_0$ be a solenoidal vector field with
coefficients in $L^{N,\infty}.$
Assume that $\theta_0$ belongs to  
\begin{itemize}
\item $L^1\cap L^{p,\infty}$ for some $p>3/2$ if $N=3,$
\item $L^{4/3,\infty}\cap L^{p,\infty}$
for some $p>4/3$ if $N=4,$
\item $L^{N/3,\infty}\cap L^{p,\infty}$ 
for some $p\geq N/3$ if $N\geq5.$
\end{itemize}
There exists a constant $c>0$ depending only on  $N,$ 
and such that if 
\begin{equation}\label{eq:petitesse}
 \|u_0\|_{L^{3,\infty}}+\|\theta_0\|_{L^{1}}\leq c\nu \ \ \hbox{if}\
\ N=3,\qquad \|u_0\|_{L^{N,\infty}}+\|\theta_0\|_{L^{N/3,\infty}}\leq c\nu \ \
\hbox{if}\ \ N\geq4 \end{equation} 
then system $\eqref{eq:boussinesq}$
has a  global solution  $(\theta,u,\nabla\Pi)$ with
$u\in L^\infty(\R_+;L^{N,\infty})$
and 
$$ \theta\in L^\infty(\R_+;L^{1}\cap L^{p,\infty})\ \ 
\hbox{if}\ \ N=3,\qquad
\theta\in L^\infty(\R_+;L^{\frac N3,\infty}\cap L^{p,\infty})
\ \ \hbox{if}\ \ N\geq4.
$$
 \end{theo} 
\remark
In fact, the heat semi-group supplies
some 
additional regularity properties for~$u$ (that we 
shall use in the proof of Theorem \ref{th:existence}). 

Note that from dimension $5$ on, one can take
$p=N/3$ so that the above statement is a global 
existence result in a scaling invariant space
for the system. 
The scaling of the system may still be almost achieved in dimension $4.$
In dimension $3$ however, it
is very unlikely that an existence result may be proved
for $\theta_0\in L^1$ and $u_0\in L^{N,\infty}$
(see Remark~\ref{r:explanation} for further explanations). 
\medbreak
 As regards uniqueness, let us stress that 
for general $u_0$ in $L^{N,\infty}$ 
the function $e^{t\Delta}u_0$
(where $(e^{t\Delta})_{t>0}$ denotes the heat semi-group)
 need not be in $L^1_{loc}(\R_+;\Lip).$
Therefore, the velocity field of the 
solution constructed in Theorem \ref{th:existence}
need not be in $L^1_{loc}(\R_+;\Lip)$
either which precludes us from 
 proving stability estimates
for system \eqref{eq:boussinesq} by mean of standard arguments.
Indeed, 
we have to deal with a transport equation 
associated to a vector field which \emph{is not}
in  $L^1_{loc}(\R_+;\Lip).$ This difficulty has been 
overcome in \cite{DANCHINPAICU} in 
the framework of  two-dimensional finite energy solutions. 
In the present paper, we shall see
that similar arguments may be used to state uniqueness
in dimension $N\geq3.$
\begin{theo}\label{th:uniqueness} 
Let $(\theta_1,u_1,\nabla\Pi_1)$ and $(\theta_2,u_2,\nabla\Pi_2)$
satisfy $\eqref{eq:boussinesq}$ with the same data. 
Assume that for some $p\in[1,2N[$ and $i=1,2,$ we have  
$$\theta_i\in L^\infty_T(\dot B^{-1+\frac 
Np}_{p,\infty})\quad\!\!\text{and}\quad\!\! 
u_i\in L_T^\infty(\dot 
B^{-1+\frac Np}_{p,\infty})\cap \tilde L_T^1(\dot B^{1+\frac Np}_{p,\infty}). 
$$ 
There exists a constant $c>0$ depending only on $N$ and on $p$ 
such that if in addition  
\begin{equation}\label{eq:petit}
\|u_1\|_{\tilde L_T^1(\dot B^{1+\frac Np}_{p,\infty})} 
+\nu^{-1}\|u_2\|_{L_T^\infty(\dot B^{-1+\frac Np}_{p,\infty})}\leq c 
\end{equation}
then  $(\theta_1,u_1,\nabla\Pi_1)\equiv(\theta_2,u_2,\nabla\Pi_2)$.
\end{theo} 
In the above statement, the space $\tilde L_T^1(\dot B^{1+\frac
Np}_{p,\infty})$ is slightly larger
than the set $L^1_T(\dot B^{1+\frac
Np}_{p,\infty})$ of integrable functions over $[0,T]$
with values in the Besov space $\dot B^{1+\frac Np}_{p,\infty}$
(see section \ref{s:tools}).
\remark
The limit case $p=2N$ may be considered if the velocity
field belongs to a Besov space with third index $1$
(see
Theorem \ref{th:uniquenesslimit}).
\medbreak
Finally, putting 
together the embedding $L^{N,\infty}\hookrightarrow
\dot B^{-1+\frac Nq}_{q,\infty}$ for $q>N$
(see Lemma \ref{l:lorentz} below),
 the existence and uniqueness theorems,
and the further regularity properties for the
velocity  given by the heat semi-group, 
one ends up with the following global well-posedness result:
\begin{theo}\label{th:GWP}  
Assume that $(\theta_0,u_0)$ satisfies
the assumptions of Theorem $\ref{th:existence}$
with $p=N.$
Then system $\eqref{eq:boussinesq}$ has
a unique solution $(\theta,u,\nabla\Pi)$
such that $u\in L^\infty(\R_+;L^{N,\infty}),$
$\theta\in L^\infty(\R_+;L^{\frac N3,\infty})$ if $N\geq4$
and $\theta\in L^1(\R_+;L^1)$ if $N=3$,
with moreover
$$
u\in\tilde L^1_{loc}(\dot B^{1+\frac Nq}_{q,\infty})\quad\hbox{for all}
\ \ q>N.
$$ 
 \end{theo}         
The paper is structured as follows.
In the next section, we present a few tools borrowed from 
harmonic and functional analysis.  
Section \ref{s:existence} is devoted
to the proof of existence. The study of uniqueness
is postponed  in section \ref{s:uniqueness}.
\medbreak\noindent{\bf Acknowledgments:}
The authors are indebted to the referee for fruitful 
remarks.

%%%%%%%%%%%%%%%%%%%%%%%%%%%%%%%%%%%%%%%%%%%%%%%%%%%%%%%%%%%%%%%%

\section{Tools and functional spaces}\label{s:tools}

\subsection{Lorentz spaces}\label{ss:lorentz}
To start with, 
let us
recall the definition of  weak $L^p$ spaces (denoted by $L^{p,\infty}$):
 \begin{defi} For $1\leq p<\infty,$ we denote by   $L^{p,\infty}(\R^N)$ (or
simply $L^{p,\infty}$) the space of all 
real valued measurable functions over
$\R^N$ such that 
$$ 
\|f\|_{L^{p,\infty}}:=\sup_{\lambda>0} \lambda\: 
\Bigl|\bigl\{x\in\R^N/|f(x)|>\lambda\bigr\}\Bigr|^{\frac1p}<\infty. 
$$ 
\end{defi} 
\begin{rem}\label{r:lorentz} 
The space $L^{p,\infty}$ 
may be alternatively defined 
by real interpolation:
$$L^{p,\infty}=(L^\infty,L^1)_{(\frac{1}{p},\infty)}.$$ 
In other words, a function 
 $f$ belongs to   $L^{p,\infty}$ if and only if, for 
all $A>0$, one may write $f=f^A+f_A$ 
for some functions $f_A\in L^1$ and $f^A$ in $L^\infty$ such that
 $$\|f_A\|_{L^1}\leq 
CA^{1-1/p}\quad\hbox{and}\quad \|f^A\|_{L^\infty}\leq C A^{-1/p}.$$
The ``best  constant" $C$ defines
a norm which is equivalent 
to $\|f\|_{L^{p,\infty}}.$ 
 \end{rem} 
The set
$\cC_c^\infty$ of smooth compactly supported 
functions is \emph{not} dense in spaces $L^{p,\infty}.$  
It turns out however  that 
$\cC_c^\infty$ is {\it  locally dense}
in  $L^{p,\infty}$ 
(despite the fact that  $L^{p,\infty}$ is not separable).
More details are given in the following proposition.
\begin{pro}\label{p:density}
 For all  $p\!\in]1,\infty[$ and $f\in L^{p,\infty},$
there exists a family $(f_\e)_{\e>0}$ of   
$\cC_c^\infty$ functions and a constant $C$ so that 
\begin{equation}\label{eq:approx}
\sup_{\e>0}\Vert f_\e\Vert_{L^{p,\infty}}\leq C\Vert f\Vert_{L^{p,\infty}}
\quad\hbox{and}\quad 
f_\e\rightarrow f\quad\hbox{in}\ \ L^1+L^\infty.
\end{equation}
\end{pro} 
\begin{p} 
Let  $(\varphi_\epsilon)_{\epsilon>0}$ be a family
of mollifiers and 
 $(\chi_\epsilon)_{\epsilon>0},$ a family
of cut-off functions with values
in $[0,1],$ supported in  $B(0,2\epsilon^{-1})$
and equal to $1$ in a neighborhood of  $B(0,\epsilon^{-1}).$

Let  $f\in L^{p,\infty}.$ For all $\e>0,$ set
$f_\e:=\varphi_\e\star(\chi_\e f)$ and fix some $\eta>0.$
{}From the above definition and remark, it is obvious that
$$
\Vert f_\e\Vert_{L^{p,\infty}}
\leq C\Vert \chi_\e f\Vert_{L^{p,\infty}}
\leq C\Vert f\Vert_{L^{p,\infty}}\quad\hbox{for all}\quad\e>0,
$$
and that one can find two functions 
$g\in L^1$ and $h\in L^\infty$ 
such that 
$f=g+h$ and $\Vert h\Vert_{L^\infty}\leq\eta/4$
(just take $A$ large enough). 

Let us split $f_\e-f$ as follows:
$$
f_\e-f=\Bigl\{\varphi_\e\star\bigl(\chi_\e g\bigr)-\varphi_\e\star g\Bigr\}
+\Bigl\{\varphi_\e\star g-g\Bigr\}
+\Bigl\{\varphi_\e\star\bigl(\chi_\e h\bigr)-h\Bigr\}.
$$
On the one hand, Lebesgue dominated convergence theorem
and standard results on convolution ensure
that the first two terms between curly brackets
have $L^1$ norm less than $\eta/2$
for small enough $\e.$
On the other hand, the $L^\infty$ norm 
of the last term is obviously less than
$\eta/2$ for any $\e>0.$

 It is now easy to complete the proof of the proposition.   \end{p} 

%%%%%%%%%%%%%%%%%%%%%%%%%%%%%%%%%%%%%%%%%%%%%%%%%

\subsection{Besov spaces}

In this subsection we  define 
 the Littlewood-Paley decomposition and 
the Besov spaces that we are going to work with.
The reader is referred to the monographs
\cite{C}, \cite{L} or \cite{RS}
for a more detailed presentation.
 \smallbreak
To start with, fix a smooth nonnegative  radial function  $\chi$ 
with support in the ball  $\{\vert\xi\vert\leq\frac{4}{3}\},$ 
value $1$ over  $\{\vert\xi\vert\leq\frac{3}{4}\},$ 
and such that  $r\mapsto \chi(re_r)$ be nonincreasing over $\R_+.$
Let  $\varphi(\xi)=\chi(\xi/2)-\chi(\xi).$
We obviously have 
\begin{equation}\label{eq:dyadique} 
\sum_{q\in \Z}\varphi(2^{-q}\xi)=1\quad\hbox{for all}\quad
\xi\in\R^N\setminus\{0\}.
\end{equation} 
We define the spectral localization operators
$\dot\Delta_q$ and  $\dot S_q$ ($q\in\Z$) by
$$ 
\dot\Delta_q\;u:=\varphi(2^{-q}D)u\quad\text{and}\quad 
\dot S_q\,u:=\chi(2^{-q}D)u. 
$$ 
For any tempered distribution 
$u\in\mathcal S'(\R^N)$, functions  $\dot\Delta_q u$ and 
$\dot S_qu$ are analytic with at most polynomial growth
and 
$\displaystyle{u=\sum_{q\in\Z}\dot\Delta_q u}$ 
 modulo  polynomials (see e.g. \cite{L}).  
\medbreak 
We shall often use the following quasi-orthogonality
property~: 
\begin{equation}\label{presortho} 
\dot\Delta_k\dot\Delta_q u\equiv 0 \quad\mbox{if}\quad\vert k-q\vert\geq 2 
\quad\mbox{and}\quad\dot\Delta_k( \dot S_{q-1}u\dot\Delta_qv)\equiv 
0\quad\mbox{if} \quad\vert k-q\vert\geq 5. 
\end{equation} 
Let us now recall the definition of
homogeneous Besov spaces. 
\begin{defi}\label{def} 
Let $s\in\R,\,(p,r)\in[1,\infty]^2$ and $u\in {\cal {S}}'(\R^N).$ 
We set
$$ 
\Vert u\Vert_{ \dot B^s_{p,r}}:= 
\Big(\sum_{q\in\Z}2^{rqs} 
\Vert\dot\Delta_q\,u\Vert_{L^p}^r\Big)^{\frac{1}{r}}\quad\text{if}\quad 
r<\infty,\quad \text{and}\!\quad \Vert u\Vert_{ \dot B^s_{p,\infty}} 
:=\sup_{q\in\Z}\,2^{qs}\Vert\dot\Delta_q\,u 
\Vert_{L^p}. 
$$ 
\begin{itemize}
\item If $s<\frac Np$ or $s=\frac Np$ and $r=1$ then 
the homogeneous  Besov space
 $\dot B^s_{p,r}:=\dot B^s_{p,r}(\R^N)$ is defined as
the set of those tempered distributions $u$ such that 
$\|u\|_{\dot B^s_{p,r}}<\infty.$
\item If $\frac{N}{p}+k\leq s<\frac{N}{p}+k+1$  (or 
$s=\frac{N}{p}+k+1$ and $r=1$) for some $k\in\N$
then   $\dot B^s_{p,r}$ is
the set of tempered distributions $u$ so that 
$\partial^\alpha u\in \dot B^{s-k-1}_{p,r}$ for 
all multi-index $\alpha$ of length $k+1.$
\end{itemize}
\end{defi} 
\remark 
Let us all also recall in passing that   
 the \emph{nonhomogeneous} Besov space
$B^s_{p,r}$ is the set of tempered distributions $u$ such 
that\footnote{note 
that the definition is somewhat simpler
than in the homogeneous framework since 
 low frequencies cannot experience divergence.}
$$
\|u\|_{B^s_{p,r}}<\infty
\quad\hbox{with}\quad
\Vert u\Vert_{B^s_{p,r}}:=
\Vert\dot S_0u\Vert_{L^p}
+\bigl\Vert 2^{qs}\Vert \ddq u\Vert_{L^p}\bigr\Vert_{\ell^r(\N)}.
 $$
\remark 
For all $s\in\R,$ the Besov space $\dot B^s_{2,2}$ coincides 
with the homogeneous  Sobolev space $\dot H^s,$ 
and, if  $s\in\R_+\setminus\N$ then  
 $\dot B^s_{\infty,\infty}$  coincides
with the homogeneous   H\"older space~$\dot C^s.$ 
\medbreak
The following \emph{Bernstein inequality}
 will be of constant use in the paper.
\begin{lem}  Let $k\in\N,$ $1\leq p_1\leq p_2\leq\infty$ and
$\psi\in\cC_c^\infty(\R^N).$ There exists a constant $C$
depending only on $k,$ $N$ and ${\rm Supp}\,\psi$ such that
$$\hskip1cm\|D^k\psi(2^{-q}D)u\|_{L^{p_2}}
\leq C2^{q\bigl(k+N\bigl(\frac 1{p_{1}}-\frac 
1{p_{2}}\bigr)\bigr)}\|\psi(2^{-q}D)u\|_{L^{p_1}}.$$ 
\end{lem} 
We shall also use the following 
fundamental properties of Besov spaces.
\begin{pro}\label{injection}
(i) The set $\dot B^s_{p,r}(\R^N)$ is a complete
 subspace of ${\cal S}'(\R^N)$ if and only if 
$s<N/p$ or $s=N/p$ and $r=1.$\\ 
(ii) There exists a positive constant  $c$ so that 
\begin{equation}\label{equivalence} 
c^{-1}\Vert u\Vert_{\dot B^s_{p,r}} \leq \Vert\nabla 
u\Vert_{\dot B^{s-1}_{p,r}} \leq c \Vert u\Vert_{\dot B^s_{p,r}}. 
\end{equation} 
(iii) For all $1\leq p_1\leq p_2\leq\infty$ and
 $1\leq r_1\leq r_2\leq\infty,$ we have 
$\dot B^s_{p_1,r_1}\hookrightarrow 
\dot B^{s-N(\frac{1}{p_1}-\frac{1}{p_2})}_{p_2,r_2}.$\\ 
(iv) If $p\in[1,\infty]$ then  
$\dot B^{\frac{N}{p}}_{p,1}\hookrightarrow 
\dot B^{\frac{N}{p}}_{p,\infty}\cap L^{\infty}.$ 
If  $p$ is finite then space $\dot B^{\frac Np}_{p,1}$ is
an algebra.\\ 
(v) Real interpolation: $(\dot B^{s_1}_{p,r}, 
\dot B^{s_2}_{p,r})_{\theta,\,r'}= \dot B^{\theta s_2+(1-\theta)s_1}_{p,r'}$ 
whenever  $0<\theta<1$ and $1\leq p,r,r'\!\leq~\!\!\infty.$ 
\end{pro}
We shall often use the fact that
Lorentz spaces
are embedded in Besov spaces (see the proof in \cite{DANCHINPAICU}).
\begin{lem}\label{l:lorentz} For any $1<p<q\leq\infty,$ we have 
$$ 
L^{p,\infty}(\R^N)\hookrightarrow \dot B^{\frac Nq-\frac 
Np}_{q,\infty}(\R^N). 
$$ 
\end{lem}
In order to pass to the limit
in the nonlinear terms of System \eqref{eq:boussinesq}, 
the following compactness result in \emph{nonhomogeneous}
Besov spaces will be most
useful.
\begin{lem}\label{l:compact}
For any $1<p<q\leq\infty,$ $\e>0$ and
$\psi\in\cC_c^\infty,$ the 
map $u\mapsto \psi u$
is compact from $L^{p,\infty}(\R^N)$
to  $B^{\frac Nq-\frac 
Np-\e}_{q,\infty}(\R^N).$
\end{lem}
\begin{p}
It is well known 
that $u\mapsto \psi u$
is a compact mapping 
from $B^s_{p,r}$ to $B^{s'}_{p,r}$
whenever $s'<s$ (see e.g
\cite{RS} and the references therein).
So it suffices to combine
the previous lemma
with the fact that, owing to $N/q-N/p<0,$ we have
$\dot B^{\frac Nq-\frac Np}_{q,\infty}
\hookrightarrow  B^{\frac Nq-\frac 
Np}_{q,\infty}.$
\end{p}

%%%%%%%%%%%%%%%%%%%%%%%%%%%%%%%%%%%%%%%%%%%%%%%%%

\subsection{Product and paraproduct in Besov spaces}
 
In order to bound the nonlinear terms in system 
\eqref{eq:boussinesq}, 
it will be useful to know how the product
of two functions operates in Besov spaces. 
In fact, optimal results may be achieved
by taking advantage of (basic) paradifferential calculus, 
a tool which was introduced by J.-M. Bony
in \cite{B}. 
More precisely, the product of two functions 
$f$ and $g$ may be decomposed according to
\begin{equation}\label{eq:decompobony}
fg=\dot T_fg+\dot T_gf+\dot R(f,g)
\end{equation}
where the paraproduct operator $\dot T$ is defined by 
the formula
$$\displaystyle{\,\dot T_fg:=\sum_q\dot S_{q\!-\!1}f\,\dot\Delta_q g},$$
and the remainder operator, $\dot R,$ by
 $$\displaystyle{\,\dot R(f,g)
:=\sum_q \dot\Delta_qf\tilde\Delta_qg}
\quad\text{with}\quad\tilde\Delta_q:=\dot\Delta_{q\!-\!1} +\dot\Delta_q
+\dot\Delta_{q\!+\!1}.$$
We shall make an extensive use of the following results
of continuity for operators $\dot T$ and $\dot R$
(see the proof in e.g. \cite{RS}):
\begin{pro}\label{p:paraproduit}
Let  $1\leq p,p_1,p_2,r,r_1,r_2\leq\infty$
so that  $\frac1p=\frac1{p_1}+\frac1{p_2}$ and
 $\frac1r=\frac1{r_1}+\frac1{r_2}.$

\noindent Operator $\dot T$ is continuous~:
\begin{itemize}
\item from $L^\infty\times\dot B^t_{p,r}$ to
$\dot B^t_{p,r}$ for all $t\in\R,$
\item from  $\dot B^{-s}_{p_1,r_1}\times
\dot B^t_{p_2,r_2}$ to $\dot B^{t-s}_{p,r}$
for all  $t\in\R$ and  $s>0.$
\end{itemize}
Operator $\dot R$ is continuous:
\begin{itemize}
\item from $\dot B^s_{p_1,r_1}\times\dot B^t_{p_2,r_2}$
to $\dot B^{s+t}_{p,r}$ for all $(s,t)\in\R^2$ such that $s+t>0,$
\item from $\dot B^s_{p_1,r_1}\times\dot B^{-s}_{p_2,r_2}$
to $\dot B^{0}_{p,\infty}$ if $s\in\R$ and
 $\frac1{r_1}+\frac1{r_2}\geq1.$
\end{itemize}
\end{pro}
Combining the above continuity results with
Bony's decomposition \eqref{eq:decompobony}, 
we get:
 \begin{pro}\label{loideproduitparticuliere} 
Let $(p,r)\in[1,\infty]^2$ and $(s_1,s_2)\in\R^2.$  
The following inequalities hold true:
\begin{itemize}
\item
If $s_1+s_2+N\inf(0,1-\frac{2}{p})>0,$ 
$s_1<\frac{N}{p}$ and 
$s_2<\frac{N}{p}$ then 
\begin{equation}\label{p_1p_2} 
\|uv\|_{\dot B^{s_1+s_2-\frac{N}{p}}_{p,r}} 
\leq C \|u\|_{\dot B^{s_1}_{p,r}}\|v\|_{\dot B^{s_2}_{p,\infty}}. 
\end{equation} 
If $s_1=\frac{N}{p}$ then 
$\|u\|_{\dot B^{\frac Np}_{p,r}}$ (resp.
$\|v\|_{\dot B^{s_2}_{p,\infty}}$)
has to be replaced with  $\|u\|_{\dot B^{\frac Np}_{p,r}\cap L^\infty}$
(resp. $\|v\|_{\dot B^{s_2}_{p,r}}$).
If $s_2=\frac Np$ then
$\|v\|_{\dot B^{\frac Np}_{p,\infty}}$ has to be replaced with 
$\|v\|_{\dot B^{\frac Np}_{p,\infty}\cap L^\infty}.$ 
\item If $s_1+s_2=0,$ 
$s_1\in(-\frac{N}{p},\frac{N}{p}]$ and $p\geq2$ then 
\begin{equation}\label{sommenulle} 
\|uv\|_{\dot B^{-\frac{N}{p}}_{p,\infty}} 
\leq C\|u\|_{\dot B^{s_1}_{p,1}}\|v\|_{\dot B^{s_2}_{p,\infty}}. 
\end{equation} 
\item 
If $p\geq 2$ then   
\begin{equation}\label{sommenulle2} 
\|uv\|_{\dot B^{-\frac{N}{p}}_{p,\infty}} 
\leq C\|u\|_{\dot
B^{-\frac Np}_{p,1}}\|v\|_{\dot B^{\frac Np}_{p,\infty}\cap L^\infty}.  \end{equation} 
 \item
If  $\vert s\vert<\frac{N}{p}$ and  $p\geq2$ or 
$-\frac{N}{p'}<s<\frac{N}{p}$ and $p<2$ then 
\begin{equation}\label{produit3} 
\Vert uv\Vert_{\dot B^s_{p,r}} \leq C\Vert u\Vert_{\dot B^s_{p,r}} \Vert 
v\Vert_ {\dot B^{\frac{N}{p}}_{p,\infty}\cap L^{\infty}}. 
\end{equation} 
 \end{itemize}
\end{pro}  

%%%%%%%%%%%%%%%%%%%%%%%%%%%%%%%%%%%%%%%%%%%%%%%%%

\subsection{Estimates for the heat and Stokes equations}

We shall often use  the following smoothing property
for the heat equation which has been stated by J.-Y. Chemin in
\cite{CH}. 
\begin{pro}\label{p:chaleur}
Let  $s\in\R,$ $1\leq p,r,\rho_1\leq\infty.$
Let $u_0\in \dot B^s_{p,r}$
and  $f\!\in\!\tilde L_T^{\rho_1}(\dot B^{s\!-\!2\!+\!\frac2{\rho_1}}_{p,r}).$
Then the heat equation
$$
\d_tu-\nu\Delta u=f,\qquad u_{|t=0}=u_0
$$
admits a unique solution $u$ in 
$\tilde L_T^\infty(\dot B^s_{p,r})\cap\tilde L_T^{\rho_1}(\dot
B^{s+\frac2{\rho_1}}_{p,r})$
and there exists a constant $C$ depending only on the
dimension $N$ so that for all  $t\in[0,T]$ and $\rho\geq\rho_1$, we have:
\begin{equation}\label{eq:estchaleur}
\nu^{\frac1\rho}\Vert u\Vert_{\tilde L_t^\rho(\dot B^{s+\frac2\rho}_{p,r})}
\leq C\biggl(\Vert u_0\Vert_{\dot B^s_{p,r}}+
\nu^{\frac1{\rho_1}-1}\Vert f\Vert_{\tilde
L_t^{\rho_1}(\dot B^{s-2+\frac{2}{\rho_1}}_{p,r})}\biggr).
\end{equation}
\end{pro}
In the above statement, spaces $\tilde L_T^\rho(\dot B^s_{p,r})$ 
are defined along the lines of  definition \ref{def} 
with 
\begin{equation}\label{eq:ltilde}
\Vert u\Vert_{\tilde L_T^\rho(\dot B^s_{p,r})}:= \Bigl\|2^{qs}\Vert
\dot\Delta_qu\Vert_{L_T^\rho(L^p)}\Bigr\|_{\ell^r(\Z)}.
\end{equation}
Note that by   virtue of Minkowski inequality, we have
\begin{equation}\label{eq:Minkowski}
\Vert u\Vert_{L_T^\rho(\dot B^s_{p,r})}\leq
 \Vert u\Vert_{\tilde L_T^\rho(\dot B^s_{p,r})}
\quad\text{for}\quad \rho\geq r,
\end{equation}
and the opposite inequality if $\rho\leq r.$
\smallbreak
We shall set 
 $\tilde L^\rho_{loc}(\R_+;\dot B^s_{p,r})=\bigcap_{T>0}
\tilde L_T^\rho(\dot B^s_{p,r}).$
\begin{rem}\label{r:tilde}
The above results of continuity for the paraproduct, remainder
and product may be easily carried out to 
$\tilde L_T^\rho(\dot B^s_{p,r})$ spaces.
The time exponents just behave according to H\"older inequality.
\end{rem}
\begin{rem}\label{r:stokes}
Since the Leray projector $\cP$ over solenoidal vector fields
maps $\dot B^s_{p,r}$ to $\dot B^s_{p,r}$ for every
$s\in\R$ and $1\leq p,r\leq\infty,$ 
Proposition $\ref{p:chaleur}$ may be extended to the nonstationary
Stokes system
\begin{equation}\label{eq:stokes}
\left\{\begin{array}{l}\d_tu-\nu\Delta u+\nabla\Pi=f,
\qquad{\rm div}\, u=0,\\
u_{|t=0}=u_0,\end{array} \right.
\end{equation}
with divergence free  initial data  $u_0\in\dot B^s_{p,r},$ 
and source term $f$ in $\tilde L_T^1(\dot B^s_{p,r}).$
In particular, denoting $\cQ={\rm Id}-\cP,$ we have the following a priori
estimates for  all $\rho\geq1:$
$$\nu^{\frac1\rho}\Vert u\Vert_{\tilde L_T^\rho(\dot B^{s+\frac2\rho}_{p,r})}
\leq C\biggl(\Vert u_0\Vert_{\dot B^s_{p,r}} +\Vert \cP f\Vert_{\tilde
L_T^1(\dot B^{s}_{p,r})}\biggr)\quad\hbox{and}\quad
\Vert\nabla\Pi\Vert_{\tilde
L_T^1(\dot B^{s}_{p,r})}\leq C\Vert\cQ f\Vert_{\tilde
L_T^1(\dot B^{s}_{p,r})}.
$$
\end{rem}

%%%%%%%%%%%%%%%%%%%%%%%%%%%%%%%%%%%%%%%%%%%%%%%%%%%%%%%%%%%%%%%

\section{Existence of solutions with infinite energy}
\label{s:existence}

This section is devoted to the proof of Theorem \ref{th:existence}. The
principle of the proof is standard:
\begin{enumerate}
\item approximate the data $(\theta_0,u_0)$
by a sequence  $(\theta_0^n,u_0^n)_{n\in\N}$ of smooth solutions;  
\item solve (globally) the Boussinesq system 
with data $(\theta_0^n,u_0^n);$ 
\item resort to compactness arguments to 
prove the convergence of a subsequence;
\item pass to the limit in the system.
\end{enumerate}

%%%%%%%%%%%%%%%%%%%%%%%%%%%%%%%%%%%%%%%%%%%%%%%%%%%%%%%%%%%%

\subsection{Global existence: the smooth case}

As a warm up, let us first  consider the case $p=N$
which is easier to deal with.
For notational simplicity, we agree that 
$L^{\frac N3,\infty}$ stands for $L^1$ if $N=3.$

Let   
$(\theta_0,u_0)$ be in  
$\bigl(L^{\frac N3,\infty}\cap L^{N,\infty}\bigr) 
\times \bigl(L^{N,\infty}\bigr)^N$ with $\div u_0=0.$
According to (a slight generalization of) Proposition
\ref{p:density}, one
can find a sequence 
$(\theta_0^n,u_0^n)_{n\in\N}$
of  ${\cal C}_c^\infty$ functions which tends 
to $(\theta_0,u_0)$ 
in the sense of the distributions 
and such that in addition  
$$ 
\|\theta_0^n\|_{L^{\frac N3,\infty}}\leq 
C\|\theta_0\|_{L^{\frac N3,\infty}},\quad 
\|\theta_0^n\|_{L^{N,\infty}}\leq 
C\|\theta_0\|_{L^{N,\infty}},\quad 
\|u_0^n\|_{L^{N,\infty}}\leq 
C\|u_0\|_{L^{N,\infty}}. 
$$ 
Note that projecting on the set of solenoidal 
vector fields is needed 
to ensure that $\div u_0^n=0.$
This operation is harmless since $\cP$ maps $L^{N,\infty}$
in $L^{N,\infty}.$
\smallbreak
According to 
  \cite{DANCHINPAICU} Theorem  1.5, we deduce
that there exists a positive constant $c$
so that if condition \eqref{eq:petitesse} is
satisfied then system \eqref{eq:boussinesq} 
admits a unique global solution $(\theta^n,u^n)$  
in  $\cC(\R_+;\dot B^0_{N,1})$ with 
in addition
$$\|u^n(t)\|_{L^{N,\infty}}\leq
 C(\nu^{-1}\|\theta_0\|_{L^{\frac N3,\infty}} 
+\|u_0\|_{L^{N,\infty}})\quad\hbox{for all}\ \ t\geq0.$$ 
As  
 $\theta_0\in L^{\frac N3,\infty}\cap L^{N,\infty}$
 and $\div u^n=0$ one can assert (see e.g Proposition 4.6 of
\cite{DANCHINPAICU}) that
$$\|\theta^n(t)\|_{L^{\frac N3,\infty}}
\leq C\|\theta_0\|_{L^{\frac N3,\infty}} 
\quad\text{and}\quad 
\|\theta^n(t)\|_{L^{N,\infty}}\leq C\|\theta_0\|_{L^{N,\infty}}
\quad\hbox{for all}\ \ t\geq0.$$ 
Note that the above estimates
are not sufficient to pass to the
limit in the nonlinear terms of system \eqref{eq:boussinesq}.
So we shall take advantage 
of the smoothing properties given by 
the Stokes  operator
in the equation for the velocity. 
Let $q$ be in $]N,\infty[.$ 
Putting together Remark~\ref{r:stokes} and the 
embedding $L^{N,\infty}\hookrightarrow\dot B^{\frac Nq-1}_{q,\infty},$ 
we discover that  
\begin{equation}\label{eq:existence1} 
\nu\|u^n\|_{\tilde L^1_t(\dot B^{\frac Nq+1}_{q,\infty})} 
\leq C\Bigl(
\|u_0\|_{L^{N,\infty}} +\|u^n\otimes u^n\|_{\tilde L^1_t(\dot B^{\frac 
Nq}_{q,\infty})} +\|\theta^n\|_{L^1_t(L^{N,\infty})}\Bigr). 
\end{equation} 
Using Bony's decomposition
\eqref{eq:decompobony} followed
by Proposition \ref{p:paraproduit}, 
Remark \ref{r:tilde}  and Lemma \ref{l:lorentz}, we can write that
$$ \begin{array}{lll}
\|u^n\otimes u^n\|_{\tilde L^1_t(\dot B^{\frac Nq}_{q,\infty})} 
&\leq& C\|u^n\|_{L_t^\infty(\dot B^{-1}_{\infty,\infty})} 
\|u^n\|_{\tilde L^1_t(\dot B^{\frac Nq+1}_{q,\infty})},\\[1.5ex] 
&\leq& C\|u^n\|_{L_t^\infty(L^{N,\infty})} 
\|u^n\|_{\tilde L^1_t(\dot B^{\frac Nq+1}_{q,\infty})}.
\end{array} 
$$ 
Plugging this latter inequality in  \eqref{eq:existence1}, 
and assuming that $c$ is suitably small, 
we thus get 
\begin{equation}\label{eq:existence2} 
\nu\|u^n\|_{\tilde L^1_t(\dot B^{\frac Nq+1}_{q,\infty})} 
\leq C\bigl(\|u_0\|_{L^{N,\infty}} +t\|\theta_0\|_{L^{N,\infty}}\bigr). 
\end{equation}  
Now, compactness arguments will enable us to pass
to the limit. Indeed, 
we have  proved that sequence 
$(\theta^n,u^n)_{n\in\N}$ is bounded in
$$ 
L^\infty(\R_+;L^{\frac N3,\infty}\cap L^{N,\infty}) 
\times\Bigl(L^\infty(\R_+;L^{N,\infty}) 
\cap\tilde L^1_{loc}(\R_+;\dot B^{\frac Nq+1}_{q,\infty})\Bigr)^N 
\quad\text{for all}\quad q>N. 
$$ 
So it is easy to show that sequence
$(\d_t\theta^n)_{n\in\N}$ is bounded in the set
of space  derivatives of functions of 
$L^\infty(\R_+;L^{\frac N2,\infty}),$ which is
embedded in  $L^\infty(\R_+;B^{-3+\frac Nq}_{q,\infty}).$ 
Now, according to Lemma \ref{l:compact}, 
for all $\psi\in\cC_c^\infty$ the map
$f\mapsto \psi f$ is \emph{compact} from    
$L^{N,\infty}$ to $B^{-3+\frac 
Nq}_{q,\infty}.$  
So one may conclude by combining
Ascoli theorem, Cantor diagonal process and interpolation that
there exists some function $\theta\in  
L^\infty(\R_+;L^{\frac N3}\cap L^{N,\infty})$ 
such that, up to  extraction,  
$$
\psi\theta^n\rightarrow\psi\theta\quad\text{in}\quad 
L^\infty_{loc}(\R_+;B^{-1+\frac Nq-\e}_{q,\infty})\quad 
\text{for all}\quad\psi\in\cC_c^\infty,
\quad\e\in]0,1[\quad\text{and}\quad q>N. 
$$ 
 Similar arguments show that
 $(u^n\otimes u^n)_{n\in\N}$ is bounded in the
space $L^\infty(\R_+;L^{\frac N2,\infty})$ 
and that   $(\theta^n)_{n\in\N}$ is bounded in 
$L^\infty(\R_+;L^{\frac N3,\infty})$ 
so that, using embeddings, we conclude that 
$(\d_tu^n)_{n\in\N}$  is bounded in   
 $L^\infty(\R_+;B^{-3+\frac Nq}_{q,\infty}).$ 
Repeating the above compactness argument, 
we get some  distribution  
$u\in L^\infty(\R_+;L^{N,\infty})\cap 
\tilde L^1_{loc}(\R_+;\dot B^{\frac Nq+1}_{q,\infty})$ 
so that, up to extraction,   
$$\psi u^n\rightarrow\psi u\quad\text{in}\quad 
L^\infty_{loc}(\R_+;B^{-1+\frac Nq-\e}_{q,\infty})\quad 
\text{for all}\quad\psi\in\cC_c^\infty,
\quad\e\in]0,1[\quad\text{and}\quad q>N. $$ 
Interpolating with the  uniform bounds in
$\tilde L^1_{loc}(\R_+;\dot B^{1+\frac Nq}_{q,\infty}),$ 
we discover that  convergence for $(\psi u^n)_{n\in\N}$ 
also holds in every space $\tilde L^r_{loc}(\R_+;B^{\frac Nq+\frac 
2r-1-\e}_{q,\infty})$ with $r>1$ and $q>N,$ 
which suffices to pass to the limit
in all the nonlinear terms.
So $(\theta, u)$ is  a (weak) solution
to  system \eqref{eq:boussinesq}. 
\hfill\rule{2.1mm}{2.1mm} 

%%%%%%%%%%%%%%%%%%%%%%%%%%%%%%%%%%%%%%%%%%%%%%%%%%%%%%%%%%%%

\subsection{Global existence: the general case}

Let us now prove the existence part
of Theorem \ref{th:existence}
in the general case. 
Let $u_0\in L^{N,\infty}$
and  
$\theta_0\in L^{\frac N3,\infty}\cap L^{p,\infty}$
for some $p$ satisfying the conditions of Theorem
\ref{th:existence}. Assume that the smallness condition \eqref{eq:petitesse}
is satisfied. 

As before, we solve  system 
\eqref{eq:boussinesq} 
with smoothed out data and  obtain 
a solution  
$(\theta^n,u^n)$ in ${\cal C}(\R_+;\dot B^0_{N,1})$ 
satisfying
$$\displaylines{\|u^n(t)\|_{L^{N,\infty}}\leq 
c(\nu^{-1}\|\theta_0\|_{L^{\frac 
N3,\infty}} +\|u_0\|_{L^{N,\infty}}),\cr 
\|\theta^n(t)\|_{L^{\frac N3,\infty}}\leq\|\theta_0\|_{L^{\frac N3,\infty}} 
\!\!\quad\text{and}\!\!\quad 
\|\theta^n(t)\|_{L^{p,\infty}}\leq\|\theta_0\|_{L^{p,\infty}}.}$$ 
Note that  $u^n$ satisfies 
$$ 
u^n(t)=\underbrace{\phantom{\!\!\!\!\!\!\!\int_0^t}e^{t\nu\Delta}u_0^n}_{u_1^n(t)}
-\underbrace{\int_0^te^{(t-\tau)\nu\Delta}{\cal P} \div(u^n\otimes
u^n)\,d\tau}_{u_2^n(t)}  +\underbrace{\int_0^t e^{(t-\tau)\nu\Delta}{\cal 
P}\bigl(\theta^n\,e_N\bigr)\,d\tau}_{u_3^n(t)}. $$ 
By virtue of the embedding
$L^{N,\infty}\hookrightarrow \dot B^{\frac N{q_1}-1}_{q_1,\infty}$ 
for  $q_1>N$ and of Proposition \ref{p:chaleur}, we get 
\begin{equation}\label{eq:u1n} 
u_1^n\in \tilde L^1_{loc}(\R_+;\dot B^{\frac N{q_1}+1}_{q_1,\infty}) 
\quad\text{uniformly with respect to}\,\ n\ \,\text{for}\ \,q_1>N. 
\end{equation} 
Next, we notice that $(u^n\otimes u^n)_{n\in\N}$ 
is bounded in  $L^\infty(\R_+;L^{\frac N2,\infty})$ 
so  that, combining the embedding  
$L^{\frac N2,\infty}\hookrightarrow \dot B^{\frac N{q_2}-2}_{q_2,\infty}$ 
for $q_2>N/2$ and  Proposition \ref{p:chaleur}, 
\begin{equation}\label{eq:u2n} 
u_2^n\!\in\!L^\infty_{loc}(\R_+;\dot B^{\frac N{q_2}-1}_{q_2,\infty}) 
\quad\text{uniformly with respect to}\ \,n\ \,\text{for}\ \, q_2>N/2. 
\end{equation} 
Finally, $(\theta^n)_{n\in\N}$ is bounded in 
 $L^\infty(\R_+;L^{p,\infty})$ 
and $L^{p,\infty}\hookrightarrow \dot B^{\frac N{q_3}-\frac 
Np}_{q_3,\infty}$ for  $q_3>p$ hence  
 \begin{equation}\label{eq:u3n} 
u_3^n\!\in\!L^\infty_{loc}(\R_+;\dot B^{\frac N{q_3}-\frac Np+2}_{q_3,\infty}) 
\quad\text{uniformly with respect to}\,\ n\ \,\text{for}\,\ q_3>p. 
\end{equation} 
Therefore, compactness arguments
similar to those which have been used in 
the previous section enable us to show that,
up to extraction, 
$(\theta^n,u^n)_{n\in\N}$ 
tends  in the sense of 
 distributions to  some $(\theta,u)$ such that 
$\theta\in L^\infty(\R_+;L^{\frac N3,\infty}\cap L^{p,\infty})$ 
and, for all $q_1>N,$ $q_2>N/2$ and $q_3>p,$ 
$$ 
u\in L^\infty(\R_+;L^{N,\infty}) 
\cap\Bigl(\tilde L^1_{loc}(\R_+;\dot B^{\frac N{q_1}+1}_{q_1,\infty}) 
+L_{loc}^\infty(\R_+;\dot B^{\frac N{q_2}-1}_{q_2,\infty}+\dot B^{\frac 
N{q_3}-\frac Np+2}_{q_3,\infty})\Bigr). $$
By interpolation with the uniform bounds
stated above, we deduce that convergence holds true
 for all  $\e\in]0,1[,$ $q>p,$ $q_1>N,$ $q_2>N/2$ and 
$q_3>p,$
\begin{itemize}
\item   
locally 
in $L^\infty_{loc}(\R_+;B^{\frac Nq-\frac Np-\e}_{q,\infty})$ 
for the temperature,
\item  locally in
$\tilde L^r\bigl(\R_+;B^{\frac
N{q_1}-1+\frac2r-\e}_{q_1,\infty}\bigr)  
+L_{loc}^\infty\Bigl(\R_+;
B^{\frac N{q_2}-1-\e}_{q_2,\infty}
+B^{\frac  N{q_3}-\frac
Np+2-\e}_{q_3,\infty}\Bigr)
 $ for the velocity.
\end{itemize}
Taking advantage of continuity 
properties for the paraproduct and the remainder,
it is then possible to pass to the limit in 
 $u^n_1\theta^n,$ $u^n_2\theta^n$ and $u^n_3\theta^n$ whenever there exist
some exponents  $q>p,$ $q_1>N,$ $q_2>N/2$ and $q_3>p$ satisfying
\begin{eqnarray}\label{eq:cond1} 
\textstyle{\frac Nq-\frac Np+\frac N{q_1}+1}&>&\textstyle{N\max\Bigl(0,\frac 1q+\frac 
1{q_1}-1\Bigr),}\\ \label{eq:cond2} 
\textstyle{\frac Nq-\frac Np+\frac N{q_2}-1}&>&\textstyle{N\max\Bigl(0,\frac
1q+\frac  1{q_2}-1\Bigr),}\\ \label{eq:cond3} 
\textstyle{\frac Nq-\frac Np+\frac N{q_3}-\frac Np+2}&>&\textstyle{N\max\Bigl(0,\frac
1q+\frac  1{q_3}-1\Bigr).} 
\end{eqnarray} 
It is clear
that the first condition is  satisfied
for $q$ close enough to $p.$
Next, second condition is verified if 
 $q_2$ is close enough to $N/2$ provided $p>N/(N-1).$ 
Finally, one can find some
$q_3$ so that the third
 condition be fulfilled if and only if 
 $p>\frac{2N}{N+2}.$ 
 
Note that in dimension three, this implies
that  $p>3/2,$ and that in dimension four, 
we must have  $p>4/3.$ 
{}From  dimension five on,
one may find some $q,$ $q_1,$ $q_2$ and $q_3$ 
such that 
 conditions \eqref{eq:cond1}, \eqref{eq:cond2}
and \eqref{eq:cond3}  are satisfied with $p=N/3.$
That no further
condition is needed to pass to the limit
in  $u^n\otimes u^n$ is left to the reader.
As a matter of fact, because $$u^n\otimes u^n=u^n\otimes(u^n_1+u^n_2+u^n_3),$$ 
it suffices to put together the fact that
 $(u^n)_{n\in\N}$ is bounded in  
$L^\infty_{loc}(\R_+;L^{N,\infty})$ 
and the properties which have been stated above for 
sequences $(u^n_1)_{n\in\N},$ 
$(u^n_2)_{n\in\N}$ and  $(u^n_3)_{n\in\N}.$ 
 \hfill\rule{2.1mm}{2.1mm} 
 \begin{rem}\label{r:explanation} 
In dimension three, one may prove
global existence under the weaker condition that $p>6/5.$ 
This may be achieved by combining the bounds satisfied 
by  $(u^n_1)_{n\in\N},$  
$(u^n_2)_{n\in\N}$ and $(u^n_3)_{n\in\N}$ 
 with a bootstrap argument. 

It is not clear however that one can take $p\leq6/5.$
Indeed, having $\theta^n$ 
in  $L^\infty(\R_+;L^{p,\infty})$ 
implies that  $u^n_3$ 
is no better than  $L^\infty_{loc}(\R_+;\dot B^{\frac 3{q_3}-\frac 
3p+2}_{q_3,\infty})$ for all $q_3>p.$ 
Hence  $p>6/5$ is needed to pass to the
 limit in  ${\rm div}\,(\theta^n\,u^n_3).$ 
\end{rem}

%%%%%%%%%%%%%%%%%%%%%%%%%%%%%%%%%%%%%%%%%%%%%%%%%%%%%%%%%%%%%

\section{Uniqueness}\label{s:uniqueness}

Let us first 
give the heuristics of the proof of
 Theorem \ref{th:uniqueness}.

Consider  two solutions 
$(\theta_1,u_1,\nabla\Pi_1)$ and $(\theta_2,u_2,\nabla\Pi_2)$
of system \eqref{eq:boussinesq} corresponding
to the same initial data. 
Assume that for some  $1\leq p<2N,$ 
$$
(\theta_i,u_i)\in L^\infty_T(\dot B^{-1+\frac Np}_{p,\infty})\times 
\Bigl(L^\infty_T(\dot B^{-1+\frac Np}_{p,\infty})\cap \tilde 
L^1_T(\dot B^{1+\frac Np}_{p,\infty})\Bigr)^N\quad\hbox{for}\ \ 
i=1,2.$$
The system  satisfied by the  difference 
$(\dt,\du,\nabla\dPi)$ between the two 
solutions reads 
$$ 
\left\{\begin{array}{l} 
\partial_t\dt+\div(u_1\dt)=-\div(\theta_2\,\du),\\ 
\partial_t\du+\div(u_1\otimes\du)-\nu\Delta\du+\nabla\dPi 
=-\div(\du\otimes u_2)+\dt\,e_N. 
\end{array}\right. 
$$ 
Note that the right-hand side of the first
equation (which is a transport equation associated
to the vector field $u_1$)
 is (at least) one derivative
less regular than $\theta_2.$
Because no smoothing property may be expected for
such an equation, this obliges us
to perform estimates in $L^\infty_T(\dot B^{-2+\frac Np}_{p,\infty})$
for $\dt$ rather than in $L^\infty_T(\dot B^{-1+\frac Np}_{p,\infty}).$
Now, due to the coupling between the equations for $\du$ and $\dt,$
this loss of one derivative also occurs in the estimates 
for $\du.$ This yields
the constraint  
$p<2N$ when bounding  
the quadratic terms (see Proposition \ref{loideproduitparticuliere}).

The second difficulty that
we have to face  is much more serious: 
as the space
$\tilde L^1_T(\dot B^{1+\frac Np}_{p,\infty})$
fails to be  embedded in 
$L^1_T(\Lip),$ the vector field $u_1$ \emph{is not}
in $L^1_T(\Lip).$ Therefore,  
the initial regularity of $\dt$ need not  be preserved
during the evolution.
It turns out however that  $\dot B^{1+\frac
  Np}_{p,\infty}$
is  embedded in the
set ${\rm Loglip}$ of  Log-Lipschitz functions
so that one may resort to 
arguments similar to those used by H. Bahouri and J.-Y. Chemin
in \cite{BCh} to prove estimates with 
(small) loss of regularity. 
Of course, we will have to cope with 
the fact that, due to the ``tilde'',
the space  $\tilde L^1_T(\dot B^{1+\frac Np}_{p,\infty})$ 
is not quite embedded in $L^1_T({\rm Loglip}).$
Overcoming this ultimate difficulty is the purpose of the next subsection.

%%%%%%%%%%%%%%%%%%%%%%%%%%%%%%%%%%%%%%%%%%%%%%%%%%%%%%%%%%%%%%%%%%%

\subsection{A priori estimates with loss of regularity}

This section is devoted to  proving
 a priori estimates  {\it with loss of regularity}
for transport-diffusion equations of the type
\begin{equation}\label{eq:td} 
\begin{cases} 
\partial_t \rho+\div(\rho u)-\nu\Delta\rho=f\\ 
\rho|_{t=0}=\rho_0 
\end{cases} 
\end{equation} 
with $u$ a given solenoidal 
vector field with coefficients
 in $\tilde L_T^1(\dot B^{1+\frac Np}_{p,\infty})$ or, 
more generally, for
the nonstationary Stokes system 
\begin{equation}\label{eq:stokesns} 
\left\{\begin{array}{l} 
\d_t v+\div(u\otimes v)-\nu\Delta v+\nabla\Pi=f,\\{\rm div}\, u=0. 
\end{array}\right. 
\end{equation}
Let us first state a result
for equation \eqref{eq:td} which is a trifle  easier to deal with. 
\begin{pro}\label{loglip}  
Let  $1<p<\infty$ and $s\in\bigl]-1-\min\bigl(\frac Np,\frac 
N{p'}\bigr),1+\frac Np\bigr[.$ Let  $\rho$ be 
a solution to the  
transport-diffusion equation
$\eqref{eq:td}$. There exists some  $N_0\in\N$ 
depending only on the choice of the Littlewood-Paley
decomposition, a universal constant  $C_0$ 
and two constants $c$ and  $C$ depending
only on  $s,$  $p$ and  $N$ so that if 
\begin{equation}\label{eq:condloglip} 
\Vert\nabla u\Vert_{\tilde L_T^1(\dot B^{\frac Np}_{p,\infty})}\leq c 
\end{equation} 
then 
 $\e_q(t):=C\Sum_{q'\leq q}2^{q'(1+\frac 
Np)}\Int_0^t\Vert\check\Delta_{q'}u\Vert_{L^p}\,d\tau$ with 
$\check\Delta_{q'}=\sum_{|\alpha|\leq N_0}\dot\Delta_{q'+\alpha}$ 
satisfies
 $$\e_q(T)-\e_{q'}(T)\leq\textstyle{\frac 12\bigl(1+s+\min\bigl(\frac 
Np,\frac N{p'}\bigr)\bigr)(q-q')}\quad\text{for all}\quad\!\! q\geq
q'$$
and  the following a priori estimate is
satisfied for all $t\in[0,T]$: 
$$\displaylines{ 
\sup_{\substack{q\in\Z\\ 
\tau\in[0,t]}}2^{qs-\e_q(\tau)}\Vert\ddq\rho(\tau)\Vert_{L^p} 
+\nu\sup_{q\in\Z}\int_0^t2^{q(s+2)-\e_q(\tau)} 
\Vert\ddq\rho(\tau)\Vert_{L^p}\,d\tau \hfill\cr\hfill\leq 
C_0\biggl(\Vert\rho_0\Vert_{\dot B^s_{p,\infty}} 
+\sup_{q\in\Z}\int_0^t2^{qs-\e_q(\tau)} \Vert\ddq 
f(\tau)\Vert_{L^p}\,d\tau \biggr).} 
$$ 
\end{pro} 
\begin{p}
The proof is similar to that  
of Proposition 4.9 
in \cite{DANCHINPAICU} (see also Theorem 3.12  
in \cite{DANCHI}).
First, we localize equation 
\eqref{eq:td} in the Fourier space
by mean of the operator $\dot\Delta_q.$ 
We get
$$ 
\d_t\ddq\rho+\dot S_{q-1}u\cdot\nabla\ddq\rho-\nu\Delta\ddq\rho=\ddq f 
+F_q
$$ 
with $F_q=F_q^1+F_q^2+F_q^3+F_q^4$ and\footnote{with the
summation convention over repeated indices $i$}
$$ 
\begin{array}{lll} 
F_q^1:=\Sum_{|q'-q|\leq4}[\dot S_{q'-1}u,\ddq]\cdot\nabla\dot\Delta_{q'}\rho,&& 
F_q^2:=\Sum_{|q'-q|\leq1}\bigl(\dot S_{q-1}-\dot S_{q'-1}\bigr)u\cdot 
\nabla\ddq\dot\Delta_{q'}\rho,\\ 
F_q^3:=-\ddq\Bigl(\Sum_{|q'\!-\!q|\leq4}\!\!\dot S_{q'-1}\d_i\rho\, 
\dot\Delta_{q'}u^i\Bigr),&& F_q^4:=-\!\Sum_{q'\geq q-3}\!\d_i\ddq\biggl( 
\dot\Delta_{q'}\rho\Bigl(\Sum_{|\alpha|\leq1}\dot\Delta_{q'\!+\!\alpha}\Bigr) 
u^i\biggr). 
\end{array} 
$$ 
Multiply both sides by $|\ddq\rho|^{p-1}\sgn(\ddq\rho),$
integrate over $\R^N$ and apply H\"older inequality. 
Owing to $\div \dot S_{q-1}u=0,$
we get:
$$
\frac 1p\frac d{dt}\Vert\ddq\rho\Vert_{L^p}^p 
-\nu\int_{\R^N}\Delta\ddq\rho\,|\ddq\rho|^{p-1}\sgn(\ddq\rho)\,dx
\leq \Bigl(\Vert\ddq 
f\Vert_{L^p} +\Vert 
F_q\Vert_{L^p}\Bigr)\Vert\ddq\rho\Vert_{L^p}^{p-1}.
$$
Now, according to Lemma A.5 in \cite{DAN}, there exists some constant
$\kappa$ depending only on $N$ and on $\Supp\varphi$
such that 
$$
-\int_{\R^N}\Delta\ddq\rho\,|\ddq\rho|^{p-1}\sgn(\ddq\rho)\,dx
\geq\kappa 2^{2q}\Vert\ddq\rho\Vert_{L^p}^p. 
$$
 Therefore, we end up with
\begin{equation}\label{eq:unicite1} 
\frac 1p\frac d{dt}\Vert\ddq\rho\Vert_{L^p}^p 
+\kappa\nu 2^{2q}\Vert\ddq\rho\Vert_{L^p}^p \leq \Bigl(\Vert\ddq 
f\Vert_{L^p} +\sum_{i=1}^4\Vert 
F_q^i\Vert_{L^p}\Bigr)\Vert\ddq\rho\Vert_{L^p}^{p-1}. 
\end{equation}
In the following calculations, assume that $2\leq p<\infty.$
Using a standard commutation estimate for bounding $F_q^1$
(see e. g. \cite{C}) and  the definition of operators $\ddq$ and $\dot S_q$
yields
$$ 
\begin{array}{ll} 
\!\!\Vert F_q^1\Vert_{L^p}\lesssim\!\Sum_{|q'\!-\!q|\leq4}\! 
\Vert{\nabla\dot S_{q'\!-\!1}u}\Vert_{L^\infty} 
\Vert\dot\Delta_{q'}\rho\Vert_{L^p},& \!\Vert 
F_q^2\Vert_{L^p}\leq\Sum_{|q'-q|\leq1} 
2^q\Vert{\check\Delta_q u}\Vert_{L^\infty}\Vert\ddq\rho\Vert_{L^p},\\ 
\!\!\Vert F_q^3\Vert_{L^p}\lesssim\!\Sum_{q'\leq q+2}\!2^{q'} 
\Vert{\dot\Delta_{q'}\rho}\Vert_{L^\infty} \Vert{\check\Delta_q 
u}\Vert_{L^p},& \!\Vert F_q^4\Vert_{L^p}\!\lesssim\!\!\Sum_{q'\!\geq\! 
q-3}\!2^{q(1\!+\!\frac Np)} 
\Vert\dot\Delta_{q'}\rho\Vert_{L^p}\Vert\check\Delta_{q'}u\Vert_{L^p}, 
\end{array} 
$$ 
with $\check\Delta_q:=\sum_{|\alpha|\leq N_0}\dot\Delta_{q'}$ 
for some large enough positive integer~$N_0.$
 \smallbreak 
Note that  
$$\displaylines{ 
\Vert{\nabla\dot S_{q'-1}u}\Vert_{L^\infty}\leq C 
\sum_{q''<q'-N_0}2^{q''(1+\frac 
Np)}\Vert\check\Delta_{q''}u\Vert_{L^p},\cr \Vert{\check\Delta_q 
u}\Vert_{L^\infty}\leq C 2^{q\frac 
Np}\Vert{\check\Delta_qu}\Vert_{L^p},\quad 
\Vert{\dot\Delta_{q'}\rho}\Vert_{L^\infty}\leq C2^{q'\frac Np} 
\Vert{\dot\Delta_{q'}\rho}\Vert_{L^p},} 
$$ 
so that plugging the above inequalities in 
\eqref{eq:unicite1}, we get
$$\displaylines{\frac 1p\frac d{dt}\Vert\ddq\rho\Vert_{L^p}^p 
+\kappa\nu2^{2q}\Vert\ddq\rho\Vert_{L^p}^p\leq\Bigl(\Vert\ddq f\Vert_{L^p} 
\hfill\cr\hfill +C\sum_{q'\geq q-4}2^{(q-q')(1+\frac Np)} 
\e'_{q'}\Vert\Delta_{q'}\rho\Vert_{L^p} +C\sum_{q'\leq 
q}2^{(q'-q)(1+\frac Np)} 
\e'_{q}\Vert\Delta_{q'}\rho\Vert_{L^p}\Bigr)
\Vert\ddq\rho\Vert_{L^p}^{p-1}} 
$$ 
with $\e_r(t):=\Int_0^t\sum_{r'\leq r}2^{r'(1+\frac Np)} 
\Vert\check\Delta_{r'}u\Vert_{L^p}\,d\tau.$ \smallbreak 
Let  $\lambda>0$ be a large enough
positive parameter (to be fixed hereafter).
We set 
$$\rho_q^\lambda(t):=2^{qs}e^{-\lambda\e_q(t)}\Vert\ddq\rho(t)\Vert_{L^p}\quad 
\text{and}\quad f_q^\lambda(t):=2^{qs}e^{-\lambda\e_q(t)}\Vert\ddq 
f(t)\Vert_{L^p}. 
$$ 
Obviously, the above inequality  rewrites 
$$\displaylines{ 
\frac 1p\frac d{dt}(\rho_q^\lambda)^p+\lambda\e'_q(\rho_q^\lambda)^p 
+\kappa\nu2^{2q}(\rho_q^\lambda)^p\leq (\rho_q^\lambda)^{p-1}
\hfill\cr\hfill
\times\biggl(f_q^\lambda+C2^{qs}e^{-\lambda\e_q} 
\Bigl(\sum_{q'\geq q-4}2^{(q-q')(1+\frac Np)} 
\e'_{q'}\Vert\Delta_{q'}\rho\Vert_{L^p} +\sum_{q'\leq 
q}2^{(q'-q)(1+\frac Np)} 
\e'_{q}\Vert\Delta_{q'}\rho\Vert_{L^p}\Bigr)\biggr)} 
$$ 
so that, performing  a time integration, we eventually get:  
$$ 
\displaylines{\rho_q^\lambda(t) 
+\kappa\nu2^{2q}\int_0^t\rho_q^\lambda(\tau)\,d\tau 
+\lambda\int_0^t\e'_q(\tau)\rho_q^\lambda(\tau)\,d\tau 
\leq\rho_q^\lambda(0)+\int_0^tf_q^\lambda(\tau)\,d\tau\hfill\cr\hfill 
+C\sum_{q'\geq q-4}2^{(q-q')(1+\frac Np+s)} 
\int_0^t\e'_{q'}(\tau)e^{\lambda(\e_{q'}(\tau)-\e_q(\tau))} 
\rho_{q'}^\lambda(\tau)\,d\tau\hfill\cr\hfill +C\sum_{q'\leq 
q}2^{(q'-q)(1+\frac Np-s)}\int_0^t 
\e'_{q}(\tau)e^{\lambda(\e_{q'}(\tau)-\e_q(\tau))} 
\rho_{q'}^\lambda(\tau)\,d\tau.} 
$$ 
Using the decomposition
 $\e'_q(\tau)=\e'_{q'}(\tau) +(\e'_q-\e'_{q'})(\tau)$ 
in the last  term and the fact that 
sequence $(\e_n)_{n\in\Z}$ is nonnegative
and nondecreasing, we gather that for all $q\in\Z,$ 
$$ 
\displaylines{\rho_q^\lambda(t)+\kappa\nu 
2^{2q}\int_0^t\rho_q^\lambda(\tau)\,d\tau 
+\lambda\int_0^t\e'_q(\tau)\rho_q^\lambda(\tau)\,d\tau 
\leq\rho_q^\lambda(0)+\int_0^tf_q^\lambda(\tau)\,d\tau 
\hfill\cr\hfill +C\sum_{q'\geq q}2^{(q-q')(1+\frac 
Np+s)}e^{\lambda(\e_{q'}(t)-\e_q(t))} 
\int_0^t\e'_{q'}(\tau)\rho_{q'}^\lambda(\tau)\,d\tau\hfill\cr\hfill 
+C\sum_{q'\leq q}2^{(q'-q)(1+\frac Np-s)}\int_0^t 
\e'_{q'}(\tau)\rho_{q'}^\lambda(\tau)\,d\tau+\frac C\lambda 
\sum_{q'\leq q}2^{(q'-q)(1+\frac Np-s)}\sup_{\tau\in[0,t]} 
\rho_{q'}^\lambda(\tau).} 
$$ 
Suppose now that the following condition is
satisfied: 
\begin{equation}\label{eq:unicite2} 
\sup_{q\in\Z}2^{q\bigl(1+\frac Np\bigr)}\Vert\check\Delta_q 
u\Vert_{L^1_t(L^p)}\leq\e\log2 \quad\text{for some }\ \e\ \text{ such 
that }\ \lambda\e\leq\frac 12\bigl(\frac Np+1+s\bigr). 
\end{equation} 
This ensures that for all $q'\geq q,$ we have 
$$ 
2^{(q-q')(1+\frac Np+s)}e^{\lambda(\e_{q'}(t)-\e_q(t))} \leq 
2^{\frac{q-q'}2(1+\frac Np+s)}. 
$$ 
Taking the supremum with respect to $q$ in 
the equation preceding \eqref{eq:unicite2},
we thus get
$$\displaylines{ 
\sup_{\substack{q\in\Z\\ 
\tau\in[0,t]}}\biggl(\rho_q^\lambda(\tau) 
+\kappa\nu2^{2q}\int_0^\tau\rho_q^\lambda(\tau')\,d\tau' 
+\lambda\int_0^\tau\e'_q(\tau')\rho_q^\lambda(\tau')\,d\tau'\biggr) 
\leq\sup_{q\in\Z}\rho_q^\lambda(0)\hfill\cr\hfill 
+\sup_{q\in\Z}\int_0^tf_q^\lambda(\tau)\,d\tau 
+C\sup_{q\in\Z}\int_0^t\e'_{q}(\tau)\rho_{q}^\lambda(\tau)\,d\tau 
+\frac C\lambda\sup_{\substack{q\in\Z\\ 
\tau\in[0,t]}}\rho_q^\lambda(\tau).} 
$$ 
Therefore,  
$$\displaylines{ 
\sup_{\substack{q\in\Z\\ 
\tau\in[0,t]}}\rho_q^\lambda(\tau) 
+\nu\sup_{q\in\Z}2^{2q}\int_0^t\rho_q^\lambda(\tau)\,d\tau 
+\lambda\sup_{q\in\Z}\int_0^t\e'_q(\tau)\rho_q^\lambda(\tau)\,d\tau 
\leq3\Vert\rho_0\Vert_{\dot B^s_{p,\infty}}\hfill\cr\hfill 
+3\sup_{q\in\Z}\int_0^tf_q^\lambda(\tau)\,d\tau 
+3C\sup_{q\in\Z}\int_0^t\e'_{q'}(\tau)\rho_{q'}^\lambda(\tau)\,d\tau 
+\frac{3C}\lambda\sup_{\substack{q\in\Z\\ 
\tau\in[0,t]}}\rho_q^\lambda(\tau).} 
$$ 
In order to conclude, it is only a matter of
choosing $\lambda=6C.$ 
We get
$$ 
\sup_{\substack{q\in\Z\\ 
\tau\in[0,t]}}\rho_q^\lambda(\tau) 
+\nu\sup_{q\in\Z}2^{2q}\int_0^t\rho_q^\lambda(\tau)\,d\tau \leq 
C_0\biggl(\Vert\rho_0\Vert_{\dot B^s_{p,\infty}} 
+\sup_{q\in\Z}\int_0^tf_q^\lambda(\tau)\,d\tau\biggr). 
$$ 
This is exactly what we wanted.
\smallbreak
In the case $p<2,$
the above bound for  $F_q^4$
turns out to be wrong. 
It may be replaced however by 
the following inequality 
$$ 
\Vert F_q^4\Vert_{L^p}\!\leq C\Sum_{q'\!\geq\!  
q-3}\!2^{q(1\!+\!\frac N{p'})}  
\Vert\Delta_{q'}\rho\Vert_{L^{p'}}\Vert\check\Delta_{q'}u\Vert_{L^p}. 
$$ 
Hence,  knowing that  
$\Vert\dot\Delta_{q'}\rho\Vert_{L^{p'}}
\leq C2^{q'(\frac Np-\frac N{p'})} 
\Vert\dot\Delta_{q'}\rho\Vert_{L^{p}},$ 
the term $2^{(q-q')(1+\frac Np)}$ 
has to replaced  by 
$2^{(q-q')(1+\frac N{p'})}$ 
in all the summations over
indices $(q,q')$ such that  $q'\geq q-4.$ 
Again, this leads to the desired inequality. 
\end{p} 
Let us now extend the previous estimate 
to the nonstationary Stokes system \eqref{eq:stokesns}.
\begin{pro}\label{loglip:stokes}  
Let $u$ be as in Proposition $\ref{loglip}$
and $v$ satisfy the nonstationary Stokes system $\eqref{eq:stokesns}.$
Then 
$$\displaylines{ 
\sup_{\substack{q\in\Z\\ 
\tau\in[0,t]}}2^{qs-\e_q(\tau)}\Vert\ddq v(\tau)\Vert_{L^p} 
+\nu\sup_{q\in\Z}\int_0^t2^{q(s+2)-\e_q(\tau)} 
\Vert\ddq v(\tau)\Vert_{L^p}\,d\tau \hfill\cr\hfill\leq 
C_0\biggl(\Vert v_0\Vert_{\dot B^s_{p,\infty}} 
+\sup_{q\in\Z}\int_0^t2^{qs-\e_q(\tau)} \Vert\ddq\cP 
f(\tau)\Vert_{L^p}\,d\tau \biggr)} 
$$ 
where $\cP$ stands for the Leray projector over solenoidal vector
fields.
\end{pro}
\begin{p}
Applying operator $\cP$ to the identity
$$ 
\d_t\ddq v+\dot S_{q-1}u\cdot\nabla\ddq v-\nu\Delta\ddq v
=\ddq f+F_q,
$$ 
and using that $\cP\ddq v=\ddq v,$ we get
$$ 
\d_t\ddq v+\dot S_{q-1}u\cdot\nabla\ddq v-\nu\Delta\ddq v
=\cP\ddq f+\cP F_q+[\dot S_{q-1}u,\cP]\cdot\nabla\ddq v.
$$ 
On the one hand,  because ${\cal F}F_q$ is supported
in some annulus $2^q C(0,r_1,r_2),$
there exists a constant $C>0$ such that
$$
\Vert{{\cal P}F_q}\Vert_{L^p}\leq C\Vert{F_q}\Vert_{L^p}.
$$
On the other hand, 
standard commutator estimates (see e.g. \cite{C}) 
ensure that  the new term $[\dot S_{q-1}u,\cP]\cdot\nabla\ddq v$
 satisfies
the same inequalities as $F_q^1.$ 
Arguing as for system \eqref{eq:td}, it is then
easy to complete the proof. 
\end{p}

%%%%%%%%%%%%%%%%%%%%%%%%%%%%%%%%%%%%%%%%%%%%%%%%%%%%%%%%%

\subsection{Proof of Theorem \ref{th:uniqueness}}\label{ss:uniqueness}

Remind that the system satisfied by the   difference 
between the two solutions reads
$$ 
\left\{\begin{array}{l} 
\partial_t\dt+\div(u_1\dt)=-\div(\theta_2\,\du),\\[1.5ex] 
\partial_t\du+\div(u_1\otimes\du)-\nu\Delta\du+\nabla\dPi 
=-\div(\du\otimes u_2)+\dt\,e_N. 
\end{array}\right. 
$$ 
We aim at proving that $(\dt,\du)\equiv0.$
To achieve it, we shall apply 
proposition \ref{loglip} or \ref{loglip:stokes} to the two equations
of the above system with  $s=-2+\frac Np-\eta$ (for
some positive $\eta$ such that  
$\eta<-1+N\min\bigl(\frac2p,1\bigr),$ 
which is consistent with the assumption $p<2N$).

So we first have to justify that 
 $\dt$ and  $\du$ belong to  
$L^\infty_T(\dot B^{-3+\frac Np+\e}_{p,\infty})$
for all $\e\in]0,1[.$
Combining
interpolation with the assumptions on the solutions
$(u_1,\theta_1)$ and $(u_2,\theta_2),$
we see that $u_i$  belongs to 
every space $\tilde L^r_T(\dot B_{p,\infty}^{\frac Np-1+\frac2r})$
with $1\leq r\leq\infty.$
Because  
$$ 
\d_t\dt=-\div\bigl(\theta_2u_2-\theta_1u_1), \qquad\dt(0)=0,
$$ 
Proposition \ref{loideproduitparticuliere} and H\"older inequality
enable us  to get  
$\dt\in{\cal C}([0,T];\dot B^{-3+\frac Np+\e}_{p,\infty})$
for all $\e\in]0,1[.$
Plugging this new information in the equation
for $\du$ and using again Proposition \ref{loideproduitparticuliere},
 it is  then easy to justify that
 $\du$ is also in $L^\infty_T(\dot B^{-3+\frac Np+\e}_{p,\infty})$.
\smallbreak  
One can now tackle the proof of 
uniqueness. Assume that the constant $c$ 
has been chosen so small as  condition 
\eqref{eq:condloglip} to be satisfied
by the vector field  $u_1.$  
Denoting
$$\displaylines{ 
\e_q(t)=C\sum_{q'\leq q}2^{q'(1+\frac 
Np)}\Vert\check\Delta_{q'}u_1\Vert_{L_t^1(L^p)},\cr 
\dTheta(t):=\sup_{\substack{\tau\in[0,t]\\q\in\Z}} 2^{-q(-2+\frac 
Np-\eta)-\e_q(\tau)} \Vert{\ddq\dt(\tau)}\Vert_{L^p},\cr 
\dU(t):=\sup_{\substack{\tau\in[0,t]\\q\in\Z}} 2^{-q(-2+\frac 
Np-\eta)-\e_q(\tau)} \Vert{\ddq\du(\tau)}\Vert_{L^p}+ 
\nu\sup_{q\in\Z}\int_0^t2^{q(\frac Np-\eta)-\e_q(\tau)} 
\Vert\ddq\du\Vert_{L^p}\,d\tau,} 
$$ 
we thus get according to Propositions \ref{loglip}
and \ref{loglip:stokes},  
$$\begin{array}{lll} 
\dTheta(t)&\!\!\!\!\leq\!\!\!\!& C\Sup_{q}\Int_0^t 2^{q(-2+\frac 
Np-\eta)-\e_q(\tau)} 
\Vert\ddq\div(\theta_2\,\du)\Vert_{L^p}\,d\tau,\\[2ex] 
\dU(t)&\!\!\!\!\leq\!\!\!\!&C\Sup_{q}\!\!\Int_0^t\!\! 
2^{q(\frac Np-2-\eta)-\e_q(\tau)} 
\bigl(\Vert\ddq\!\div(\du\otimes u_2)\Vert_{L^p} 
\!+\!\Vert\ddq\dt\Vert_{L^p}\bigr)\,d\tau.\end{array}$$ 
Let us admit that the nonlinear 
terms may be bounded as follows
(see the proof in the appendix):
$$ 
\displaylines{ 
\Sup_{q}\Int_0^t 2^{q(-2+\frac Np-\eta)-\e_q(\tau)} 
\Vert\ddq\div\!(\theta_2\,\du)\Vert_{L^p}\,d\tau\hfill\cr\hfill\leq C 
\Vert\theta_2\Vert_{L_t^\infty(\dot B^{-1+\frac Np}_{p,\infty})} 
\Sup_q\Int_0^t 2^{q(\frac Np-\eta)-\e_q(\tau)} 
\Vert\ddq\du\Vert_{L^p}\,d\tau,} 
$$ 
$$ 
\displaylines{ 
\Sup_{q}\Int_0^t\! 2^{q(-2+\frac Np-\eta)-\e_q(\tau)} 
\Vert\ddq\div\!(\du\!\otimes\! u_2)\Vert_{L^p} 
\,d\tau\hfill\cr\hfill\leq C\Vert u_2\Vert_{L_t^\infty(\dot B^{\frac 
Np-1}_{p,\infty})} \Sup_{q}\int_0^t2^{q(\frac 
Np-\eta)-\e_q(\tau)} \Vert\ddq\du\Vert_{L^p}\,d\tau.} 
$$ 
We eventually get
$$\begin{array}{lll} 
\dTheta(t)&\leq& C\nu^{-1} 
\Vert\theta_2\Vert_{L_t^\infty(\dot B^{-1+\frac Np}_{p,\infty})} 
\dU(t),\\[1ex] 
\dU(t)&\leq& C\Bigl(\nu^{-1} 
\Vert u_2\Vert_{L_t^\infty(\dot B^{\frac 
Np-1}_{p,\infty})} \dU(t)+\Int_0^t\dTheta(\tau)\,d\tau\Bigr),
\end{array} 
$$ 
whence, if $c$ has been chosen small enough in \eqref{eq:petit},   
$$ 
\dU(t)\leq C\nu^{-1}\int_0^t
\Vert\theta_2\Vert_{L_\tau^\infty(\dot B^{\frac
    Np-1}_{p,\infty})}\dU(\tau)\,d\tau. 
$$ Gronwall
Lemma  thus ensures that
$\dU\equiv0$ on $[0,T].$ Of course, this
also entails that $\dt\equiv0$ on $[0,T].$\hfill\rule{2.1mm}{2.1mm}  

%%%%%%%%%%%%%%%%%%%%%%%%%%%%%%%%%%%%%%%%%%%%%%%%%%%%%%%%%%%%%%% 
 
\subsection{The limit case $p=2N$}  

Carrying out the method which has been used in 
the previous section to the limit
case $p=2N$ seems hopeless. 
Indeed, 
 we would have to deal  with the
product $\theta_2\du$
while  $\theta_2\in L^\infty_T(\dot B^{-\frac12-\e}_{2N,\infty})$
and $\du\in \tilde L^1_T(\dot B^{\frac 12-\e}_{2N,\infty})$
for some $\e>0,$ which does not make sense
since the sum of indices of regularity is negative.

In order to meet the index $p=2N,$ one may 
resort to 
Besov spaces
\emph{with third  index~$1$}  so as to have
a velocity field 
in $L^1_T(\Lip).$ 
Thus no losing a priori estimate is needed.
On the other hand, due to the weak regularity
assumptions, 
we shall be in the limit case for product laws
in Proposition \ref{loideproduitparticuliere}
so that a  logarithmic interpolation
argument (similar to that which has been used in \cite{D} and
\cite{DANCHIN}) will be required.  Let us state the result. 
\begin{theo}\label{th:uniquenesslimit}
Let $(\theta_1,u_1,\nabla\Pi_1)$ and  $(\theta_2,u_2,\nabla\Pi_2)$ 
satisfy $\eqref{eq:boussinesq}$ with the same data.
Assume that 
$$\theta_i\in L^\infty_T(\dot B^{-\frac 1 2}_{2N,\infty})\quad\!\!\text{and}\quad\!\! 
u_i\in L_T^\infty(\dot 
B^{-\frac 12}_{2N, 1})\cap  L_T^1(\dot B^{\frac 32}_{2N,1}) 
\quad\text{for}\ \, i=1,2. 
$$ 
Then $(\theta_1,u_1,\nabla\Pi_1)\equiv(\theta_2,u_2,\nabla\Pi_2)$
on $[0,T].$ 
\end{theo} 
\begin{p}
We omit the proof of the fact that
$(\dt,\du,\nabla\dPi)$ belongs to the space
$$ 
G_T:=L^\infty_T(\dot B^{-\frac{3}{2}}_{2N,\infty}\big) 
\times\Big(\cC\big([0,T];\,\dot B^{-\frac{3}{2}}_{2N,\infty}\big)\cap 
\tilde L^1_T\big(\dot B^{\frac{1}{2}}_{2N,\infty}\big)\Big)^N 
\times\Big(\tilde L^1_T\big(\dot B^{-\frac{3}{2}}_{2N,\infty}\big)\Big)^N. 
$$ 
In the following computations, 
the space $G_T$ will be endowed with the norm
$$ 
\Vert(\theta,u,\nabla\Pi)\Vert_{G_T}\overset{d\acute{e}f}{=} \Vert 
\theta\Vert_{L^\infty_T(\dot B^{-\frac{3}{2}}_{2N,\infty})} +\Vert 
u\Vert_{L^\infty_T(\dot B^{-\frac{3}{2}}_{2N,\infty})} +\Vert 
u\Vert_{\tilde L^1_T(\dot B^{\frac{1}{2}}_{2N,\infty})} 
+\Vert\nabla\Pi\Vert_{\tilde L^1_T(\dot B^{-\frac{3}{2}}_{2N,\infty})}. 
$$ 
A  priori estimates for the transport
equation (see e.g. the limit case in Proposition 4.7 of
\cite{DANCHINPAICU}) and inequality 
\eqref{sommenulle} guarantee that for all
$t\leq T$ 
\begin{equation}\label{densite} 
\Vert\dt\Vert_{L^\infty_t(\dot B^{-\frac{3}{2}}_{2N,\infty})} \leq 
C\exp\Big(C\Vert\nabla 
u_1\Vert_{L^1_t(\dot B^{\frac{1}{2}}_{2N,1})}\Big)
\int_0^t\Vert\du\Vert_{\dot B^{\frac{1}{2}}_{2N,1}} \Vert 
\theta_2\Vert_{\dot B^{-\frac{1}{2}}_{2N,\infty}}\,d\tau. 
\end{equation} 
Next, we have according to  Proposition 3.2 of \cite{AP}
and Inequality \eqref{eq:Minkowski},  
$$\displaylines{ 
\Vert\du\Vert_{L^\infty_t(\dot B^{-\frac{3}{2}}_{2N,\infty})} 
+\nu\Vert\du\Vert_{\tilde L^1_t(\dot B^{\frac{1}{2}}_{2N,\infty})} 
+\Vert\nabla\dPi\Vert_{\tilde L^1_t(\dot B^{-\frac{3}{2}}_{2N,\infty})} 
\hfill\cr\hfill 
\leq C e^{C\Vert u_1\Vert_{L^1_t(\dot B^{\frac{3}{2}}_{2N,1})}} 
\int_0^t\Vert\dt e_N-\div(\du\otimes u_2) 
\Vert_{\dot B^{-\frac{3}{2}}_{2N,\infty}}\,d\tau.} 
$$ 
Because  $\div\du=0,$ we have for $1\leq i\leq N,$
$$
\bigl(\div(\du\otimes u_2)\bigr)^i
=\d_j\bigl(\dot T_{\du^j}u_2^i
+\dot R(\du^j,u_2^i)\bigr)+\dot T_{\d_ju_2^i}\du^j.
$$
So using Proposition \ref{p:paraproduit}, 
we find that
$$ 
\Vert\div(\du\otimes u_2)\Vert_{\dot   B^{-\frac{3}{2}}_{2N,\infty}}
\leq C\Vert\du\Vert_{\dot B^{-\frac{3}{2}}_{2N,\infty}} 
\Vert u_2\Vert_{\dot B^{\frac{3}{2}}_{2N,1}}.
$$ 
Now, taking advantage of Proposition 1.8 in \cite{D},
one may write 
$$ 
\|\du\|_{L^1_t(\dot B^{\frac{1}{2}}_{2N,1})} \leq C\|\delta 
u\|_{\tilde L^1_t(\dot B^{\frac{1}{2}}_{2N,\infty})} \log 
\bigg(e+\frac{\|\du\|_{L^1_t(\dot B^{-\frac{1}{2}}_{2N,1})} 
 +\|\du\|_{L^1_t(\dot B^{\frac{3}{2}}_{2N,1})}} 
 {\|\du\|_{\tilde L^1_t(\dot B^{\frac{1}{2}}_{2N,\infty})}}\bigg). 
$$  
Let us introduce the notation 
$$ 
\begin{array}{lll} 
V(t)&=& t\big(\|u_1\|_{L^\infty_t(\dot B^{-\frac{1}{2}}_{2N,1})} 
+\|u_2\|_{L^\infty_t(\dot B^{-\frac{1}{2}}_{2N,1})}\big) 
+\|u_1\|_{L^1_t(\dot B^{\frac{3}{2}}_{2N,1})}+ 
 \|u_2\|_{L^1_t(\dot B^{\frac{3}{2}}_{2N,1})},\\ 
W(t)&=&\Vert\du\Vert_{L^\infty_t(B^{-\frac{3}{2}}_{2N,\infty})} 
+\Vert\du\Vert_{\tilde L^1_t(\dot B^{\frac{1}{2}}_{2N,\infty})} 
+\Vert\nabla\delta\Pi\Vert_{\tilde L^1_t(\dot 
  B^{-\frac{3}{2}}_{2N,\infty})}. 
\end{array} 
$$ 
Because  
$$ 
\|\du\|_{L^1_t(\dot B^{-\frac{1}{2}}_{2N,1})}  
+\|\du\|_{L^1_t(\dot B^{\frac{3}{2}}_{2N,1})} \leq V(t), 
$$ 
and the map $x\longmapsto x\ln(e+\frac{y}{x})$
(for fixed $y\geq0$) is nondecreasing over $\R_+,$
we end up with
$$ 
\displaylines{ 
 W(t) 
\leq Ce^{C\Vert u_1\Vert_{L^1_t(\dot B^{\frac{3}{2}}_{2N,1})}} 
\int_0^t\Bigl(\|u_2(s)\|_{\dot B^{\frac32}_{2N,1}}+\|\theta_2(s)\|_{\dot
B^{-\frac12}_{2N,\infty}}\Bigr) 
W(s)\log\bigg(e+\frac{V(s)}{W(s)}\bigg)\,ds.}$$ 
Applying Osgood lemma  (see e.g. \cite{C})
thus yields $W\equiv0$ on $[0,T]$
whence also  
 $\dt=0$ according to inequality (\ref{densite}).
This completes the proof. 
 \end{p}

%%%%%%%%%%%%%%%%%%%%%%%%%%%%%%%%%%%%%%%%%%%%%%%%%%%%%%%%%%%%%
 
\section{Appendix}

This appendix is devoted to proving the  
estimates for the convection terms
that we used in section \ref{ss:uniqueness}.
\begin{lem}\label{l:quad} 
Let  $(\alpha_q)_{q\in\Z}$ be 
a sequence of nonnegative 
functions over $[0,T].$ Let  $s_1,$ $s_2,$ $p$ satisfy 
 $$1\leq p\leq\infty,\quad 
\frac Np+1>s_1,\quad 
\frac Np>s_2\ \ \text{and}\ \ s_1+s_2>N\max\Bigl(0,\frac2p-1\Bigr). 
$$ Assume that for all  $q'\geq q$ and 
 $t\in[0,T],$ we have
\begin{equation}\label{eq:alpha} 
0\leq\alpha_{q'}(t)-\alpha_{q}(t)\leq\frac12\biggl(s_1+s_2 
+N\min\Bigl(0,1-\frac2p\Bigr)\biggr)(q'-q). 
\end{equation} 
Then for all  $r\in[1,\infty],$ 
there exists a  constant $C$ 
depending only on $s_1,$ $s_2,$ $N$ and~$p$ such that for
all function $b$ and solenoidal vector field
$a$ over $\R^N,$ the following estimate  holds true
for all $t\in[0,T]$: 
$$ 
\sup_{q\in\Z} \int_0^t2^{q(s_1+s_2-1-\frac Np)-\alpha_q(\tau)} 
\Vert\ddq{\rm div}(ab)\Vert_{L^p}\,d\tau \leq C\Vert b\Vert_{\tilde 
L_t^r(\dot B^{s_1}_{p,\infty})} \sup_{q\in\Z}\Bigl\Vert 2^{q 
s_2-\alpha_q} \Vert\ddq a\Vert_{L^p}\Bigr\Vert_{L_t^{r'}}. 
$$ 
\end{lem} 
\begin{p} 
The proof relies on Bony's  decomposition.
Knowing that  $\div b=0,$ we have  (with the usual 
summation convention over repeated indices): 
\begin{equation}\label{eq:bony} 
\ddq\div(ab)=\ddq\bigl(\dot T_{\d_jb}a^j\bigr) +\ddq\bigl(\dot T_{a^j}\d_jb\bigr) 
+\ddq\d_j\dot R(a^j,b). 
\end{equation} 
By virtue of \eqref{presortho}, one may write
$$ 
\ddq\bigl(\dot T_{\d_jb}a^j\bigr)=\sum_{|q'-q|\leq4} 
\ddq\bigl(\dot S_{q'-1}\d_jb\,\dot\Delta_{q'}a^j\bigr). 
$$ 
For the sake of simplicity, let us proceed as
if 
$\ddq\bigl(T_{\d_jb}a^j\bigr)=\dot S_{q-1}\d_jb\,\ddq a^j$ 
(having \eqref{eq:alpha} justifies this
 approximation).
Using the definition of $\dot S_{q-1},$ we get 
$$ 2^{q\bigl(s_1-1-\frac Np\bigr)}
\Vert S_{q-1}\d_jb\,\ddq a^j\Vert_{L^p}\leq 
\Vert\ddq a\Vert_{L^p} \sum_{q'\leq 
  q-2}2^{q'(s_1\!-\!1\!-\!\frac Np)}\Vert\dot\Delta_{q'}\nabla 
b\Vert_{L^\infty} 2^{(q-q')(s_1\!-\!1\!-\!\frac Np)}. 
$$ 
In consequence,  we have for all $0\leq t\leq T,$   
$$\displaylines{ 
\int_0^t2^{q(s_1+s_2-1-\frac Np)-\alpha_q(\tau)} \Vert 
\dot S_{q-1}\d_jb\,\ddq a^j\Vert_{L^p}\,d\tau 
\hfill\cr\hfill\leq\Bigl\Vert 2^{qs_2-\alpha_q} \Vert\ddq 
a\Vert_{L^p}\Bigr\Vert_{L_t^{r'}} \sum_{q'\leq q-2}2^{q'(s_1-1-\frac 
Np)} \Vert\dot\Delta_{q'}\nabla b\Vert_{L_t^r(L^\infty)} 
2^{(q-q')(s_1-1-\frac Np)}.} 
$$ 
Since  $s_1-1-\frac Np<0,$ we thus get 
\begin{equation}\label{eq:Tba}\begin{array}[b]{r} 
\!\!\!\displaystyle{\sup_{q\in\Z} \int_0^t2^{q(s_1+s_2-1-\frac 
Np)-\alpha_q(\tau)} \Vert \ddq\dot T_{\d_jb}a^j\Vert_{L^p}\,d\tau} 
\\[1.5ex] 
\!\!\!\leq C\Vert \nabla b\Vert_{\tilde L_t^r(\dot B^{s_1-1-\frac 
Np}_{\infty,\infty})} 
\end{array}\sup_{q\in\Z}\Bigl\Vert 2^{qs_2-\alpha_q} 
\Vert\ddq a\Vert_{L^p}\Bigr\Vert_{L_t^{r'}}. 
\end{equation} 
Likewise, in order to bound the second term of \eqref{eq:bony}, one 
may proceed as if 
$$\ddq\bigl(\dot T_{a^j}\d_jb\bigr)=\dot S_{q-1}a^j\,\ddq\d_jb.$$ 
 Now, for all $0\leq\tau\leq T,$ we have
$$\displaylines{ 
2^{q(s_1+s_2-1-\frac Np)-\alpha_q(\tau)} \Vert 
\dot S_{q-1}a^j\ddq\d_jb\Vert_{L^p}\hfill\cr\hfill\leq 
2^{q(s_1-1)}\Vert\ddq\nabla b\Vert_{L^p} \sum_{q'\leq 
q-2}2^{-q'(\frac Np-s_2)-\alpha_q(\tau)} 
\Vert\dot\Delta_{q'}a\Vert_{L^\infty} 2^{(q'-q)(\frac Np-s_2)}.} 
$$ 
Knowing that  $\alpha_{q}\geq\alpha_{q'}$ for $q\geq q',$ 
we deduce that  
$$\displaylines{ 
\int_0^t2^{q(s_1+s_2-1-\frac Np)-\alpha_q(\tau)} \Vert 
\dot S_{q-1}a^j\ddq\d_jb\Vert_{L^p}\,d\tau 
\hfill\cr\hfill\leq2^{q(s_1-1)}\Vert\ddq\nabla b\Vert_{L_t^r(L^p)} 
\sum_{q'\leq q-2} \Bigl\Vert 2^{-q'(\frac Np-s_2)-\alpha_{q'}} 
\Vert\dot\Delta_{q'} a\Vert_{L^\infty}\Bigr\Vert_{L_t^{r'}} 
2^{(q'-q)(\frac Np-s_2)}.} 
$$ 
Combining Bernstein inequality and the fact that $s_2<\frac Np,$ 
we thus conclude that 
\begin{equation}\label{eq:Tab}\begin{array}[b]{r} 
\!\!\!\displaystyle{\sup_{q\in\Z} \int_0^t2^{q(s_1+s_2-1-\frac Np) 
-\alpha_q(\tau)} \Vert 
S_{q-1}a^j\ddq\d_jb\Vert_{L^p}\,d\tau\!\!\!} 
\\[1.5ex]\leq 
C\Vert\nabla b\Vert_{\tilde L_t^r(B^{s_1-1}_{p,\infty})} 
\end{array}\sup_{q\in\Z}\Bigl\Vert 2^{qs_2-\alpha_q} 
\Vert\ddq a\Vert_{L^p}\Bigr\Vert_{L_t^{r'}}. 
\end{equation} 
In order to treat the remainder term, we shall 
consider the cases  $p\geq2$ and $p<2.$
Let us start with the case $p\geq2$ 
which is slightly easier. We have 
$$ 
\ddq\d_j\dot R(a^j,b)=\sum_{q'\geq 
  q-3}\d_j\ddq\bigl(\dot\Delta_{q'}a^j 
\tilde\Delta_{q'}b)\quad\text{with}\quad 
\tilde\Delta_{q'}=\dot\Delta_{q'-1}+\dot\Delta_{q'}+\dot\Delta_{q'+1}.
$$ 
Because ${\cal F}\bigl(\d_j\ddq\bigl(\Delta_{q'}a^j\tilde 
\Delta_{q'}b)\bigr)$ is supported
in a ball of size $2^q,$   Bernstein
inequality ensures that 
$$ 
2^{q(s_1+s_2-1-\frac Np)-\alpha_q(\tau)}\Vert\ddq\d_j 
\dot R(a^j,b)\Vert_{L^p} \leq C\sum_{q'\geq q-3} 
2^{q(s_1+s_2)-\alpha_q(\tau)}\Vert\dot\Delta_{q'}a\Vert_{L^p} 
\Vert\tilde\Delta_{q'}b\Vert_{L^p}, 
$$ 
whence
$$\displaylines{ 
\int_0^t2^{q(s_1+s_2-1-\frac Np)-\alpha_q(\tau)}\Vert\ddq\d_j 
\dot R(a^j,b)\Vert_{L^p}\,d\tau \hfill\cr\hfill
\leq C\sum_{q'\geq 
q-3}\int_0^t 
2^{q's_2-\alpha_{q'}(\tau)}\Vert\Delta_{q'}a\Vert_{L^p}\, 
\bigl(2^{q's_1}\!\Vert\tilde\Delta_{q'}b\Vert_{L^p}\bigr)\, 
2^{(\alpha_{q'}\!-\!\alpha_q)(\tau)}2^{(q-q')(s_1+s_2)}.} 
$$ 
Thanks to assumption \eqref{eq:alpha}, we have
$$
(\alpha_{q'}-\alpha_q)(\tau)+(q-q')(s_1+s_2)\leq
(q-q')\Bigl(\frac{s_1+s_2}2\Bigr)\quad\hbox{for all}\quad
q'\geq q.
$$
As $s_1+s_2>0,$  we end up for all $q\in\Z$
with  
$$ 
\int_0^t2^{q(s_1\!+\!s_2\!-\!1\!-\!\frac Np)-\alpha_q(\tau)}\Vert\ddq\d_j 
\dot R(a^j,b)\Vert_{L^p}\,d\tau \lesssim \Vert b\Vert_{\tilde 
L_t^r(\dot B^{s_1}_{p,\infty})} \sup_{q\in\Z}\Bigl\Vert 
2^{qs_2-\alpha_q} \Vert\ddq a\Vert_{L^p}\Bigr\Vert_{L_t^{r'}}, 
$$ 
which, together with \eqref{eq:bony}, \eqref{eq:Tba} and \eqref{eq:Tab}, 
completes the proof in the case $p\geq2.$ 
\smallbreak 
If $p<2,$ one may write (use Bernstein inequality)
$$  
2^{q(s_1\!+\!s_2\!-\!1\!-\!\frac Np)-\alpha_q(\tau)}\Vert\ddq\d_j  
\dot R(a^j,b)\Vert_{L^p} \leq C\!\sum_{q'\geq q-3}  
2^{q(s_1\!+\!s_2\!+\! 
\frac N{p'}\!-\!\frac{N}p)-\alpha_q(\tau)}\Vert\dot\Delta_{q'}a\Vert_{L^{p'}}  
\Vert\tilde\Delta_{q'}b\Vert_{L^p}  $$  
and use that   
$\Vert\dot\Delta_{q'}a\Vert_{L^{p'}}\leq C2^{q'\bigl(\frac N{p}-\frac 
N{p'}\bigr)} \Vert\dot\Delta_{q'}a\Vert_{L^{p}}.$ 
Therefore, one may go along the lines
of the case $p\geq2.$ It is only a matter of changing   
the term  $2^{(q-q')(s_1+s_2)}$ 
 into $2^{(q-q')(s_1+s_2+\frac N{p'}-\frac Np)}.$ 
\end{p}

\end{document}